\documentclass[12pt]{article}
\usepackage{color}
\usepackage{amsmath,amsthm,amsfonts,amssymb}

\newtheorem{hypA}{Hypothesis}

\usepackage{mathtools}

\usepackage[hidelinks]{hyperref}
\numberwithin{equation}{section}

\newtheorem{theorem}{Theorem}[section]
\newtheorem{corollary}[theorem]{Corollary}

\newtheorem{proposition}[theorem]{Proposition}
\newtheorem{lemma}[theorem]{Lemma}

\newtheorem{remark}[theorem]{Remark}

\def\cD{\mathcal{D}}

\def\cF{\mathcal{F}}

\def\cH{\mathcal{H}}

\def\cK{\mathcal{K}}

\def\cM{\mathcal{M}}

\def\cS{\mathcal{S}}

\def\bC{\mathbb{C}}
\def\bD{\mathbb{D}}
\def\bE{\mathbb{E}}

\def\bN{\mathbb{N}}
\def\bP{\mathbb{P}}
\def\bR{\mathbb{R}}

\newcommand{\ud}{\ensuremath{\mathrm{d}}}

\newcommand{\fg}{\mathfrak{g}}
\newcommand{\fm}{\mathfrak{m}}

\newcommand{\fH}{\mathfrak{H}}
\newcommand{\fX}{\mathfrak{X}}

\begin{document}

\title{Exact super-exponential moment growth for the one-dimensional hyperbolic Anderson model: Skorohod regime}

\author{Raluca M. Balan\footnote{Corresponding author. University of Ottawa, Department of Mathematics
and Statistics, 150 Louis Pasteur Private, Ottawa, Ontario, K1N 6N5, Canada.
E-mail address: rbalan@uottawa.ca.} \and Le Chen \footnote{Auburn University,
Department of Mathematics and Statistics, Office 3147, 305 W. Samford Avenue,
Auburn, Alabama 36849, USA. E-mail address: le.chen@auburn.edu} \and Jinxin
Wang\footnote{University of Ottawa, Department of Mathematics and Statistics,
150 Louis Pasteur Private, Ottawa, Ontario, K1N 6N5, Canada. E-mail address:
jwang023@uottawa.ca.} }

\maketitle

\begin{abstract}
  \noindent In this article, we study the hyperbolic Anderson model in dimension $d=1$, driven by a Gaussian noise which is colored in space and time, with covariance $Q(t,x)=|t|^{-\beta}|x|^{-\alpha}$ for some $\alpha,\beta \in (0,1)$. The solution is interpreted in the Skorohod sense. We establish the exact moment asymptotics of the $L^p(\Omega)$-norm $\|u(t,x)\|_p$ of the solution $u$. More precisely, we show that $\|u(t,x)\|_p$ grows super-exponentially at the speed $t_p^{(4-\alpha-\beta)/(3-\alpha)}$, when $t_p=(p-1)^{1/(4-\alpha-\beta)} t$ goes to infinity. For this, we use a connection with a deterministic wave equation and a Cameron-Martin space.
\end{abstract}


\vspace{1mm}


\tableofcontents

\section{Introduction}

In this article, we consider the hyperbolic Anderson model in dimension $d = 1$:
\begin{align}
  \label{E:HAM}
  \begin{cases}
    \dfrac{\partial^2 u}{\partial t^2} (t,x)
    = \dfrac{\partial^2 u}{\partial x^2} (t,x)+\sqrt{\theta} u(t,x)
    \dot{W}(t,x), \quad
    t>0, \ x \in \bR, \\[1em]
    u(0,x) = 1, \quad \dfrac{\partial u}{\partial t} (0,x)=0, \quad x \in \bR,
  \end{cases}
\end{align}
where $\theta > 0$ and $\{W(\varphi);\varphi \in \cD\}$ (with $\cD \coloneqq
C_c^{\infty}((0,\infty) \times \bR)$) is a centered Gaussian process defined on
a complete probability space $(\Omega,\cF,\bP)$, with covariance:
\begin{equation}\label{E:cov-W}
  \bE[W(\varphi)W(\psi)]=\int_{((0,\infty) \times \bR)^2}
  \varphi(t,x)\psi(s,y)|t-s|^{-\beta}|x-y|^{-\alpha} \ud t\ud x\ud s
  \ud y=:\langle \varphi,\psi \rangle_{\cH},
\end{equation}
for some $\alpha,\beta \in (0,1)$. %

By analogy with d'Alembert's formula, we say that a process
$\{u_{\theta}(t,x);(t,x) \in [0,\infty) \times \bR\}$ is a (Skorohod) {\em
solution} of equation \eqref{E:HAM} if for any $t \geq 0$ and $x \in \bR$,
$u_\theta(t,x)$ is $(\cF_t)$-measurable, and
\[
  u_{\theta}(t,x)=1+\sqrt{\theta}\int_0^t \int_{\bR} G(t-s,x-y)
  u_{\theta}(s,y)W(\delta s,\delta y) \quad \mbox{a.s.}
\]
Here, we let $(\cF_t)_{t \geq 0}$ be the filtration induced by $W$, we interpret
the stochastic integral as a Skorohod integral (defined in
Section~\ref{S:section-background} below), and we let $G$ be the fundamental
solution of the wave equation on $(0,\infty) \times \bR$:
\[
  G(t,x)=\frac{1}{2}1_{\{|x|<t\}}.
\]

From~\cite{balan:12:stochastic}, it is known that equation~\eqref{E:HAM} has a
unique solution $u_{\theta}$. The goal of this paper is to derive the exact
asymptotics of $\bE|u_\theta(t,x)|^p$, when $t$ (or $p$) is large. This problem
has been investigated for more than a decade by the first two authors, Jian
Song, Xia Chen, and other researchers. Despite sustained efforts over the years,
the problem has remained open until now. The paper now solves this problem with
the help of AI.

\medskip

We begin by recalling some known results from the literature. We denote by
$\|\cdot\|_p$ the norm in $L^p(\Omega)$, for any $p \geq 1$.

\medskip

In~\cite{dalang.mueller:09:intermittency}, Dalang and Mueller considered the
analogue of equation \eqref{E:HAM} in dimension $d=3$ (in which
$\frac{\partial^2}{\partial x^2}$ is replaced by the Laplacian $\Delta = \sum_{i
= 1}^d \frac{\partial^2}{\partial x_i^2}$), driven by a Gaussian noise which is
white in time and has bounded spatial covariance $f$. For positive constant
initial displacement $u_0$ and velocity $\widetilde u_0$, their upper bound for
even moments gives, for $t\ge0$, $x\in\bR^3$ and even integers $p\ge2$,
\[
  \|u(t,x)\|_p \le c_1 (u_0+t\widetilde u_0)
  \exp\big(c_2\,p^{1/3}t\big).
\]
They also obtained a lower bound of the same order under a local positivity
assumption on $f$. For both bounds, the authors used a Feynman-Kac (FK)-type
representation for the moments of the solution, which they had developed in an
earlier paper with Tribe (see~\cite{dalang.mueller.ea:08:feynman-kac-type}),
using a classical Poisson process (independent of the noise). Other results for
the moments of the solution of the stochastic wave equation in any dimension $d
\geq 1$ can be found in~\cite{conus.dalang:08:non-linear, balan.song:19:second}
for the white noise in time, and in~\cite{balan.song:17:hyperbolic} for the
colored noise in time.

It turns out that a simpler method exists for the upper bound, which uses the
Wiener chaos expansion: the proof of Theorem 4.1 of~\cite{balan:12:stochastic}
shows that if $u_{\theta}$ solves the analogue of equation \eqref{E:HAM} in
dimension $d\ge1$, with spatial covariance $|x|^{-\alpha}$,
$0<\alpha<\min(d,2)$, and temporal covariance
$|t|^{-\beta}$, $0<\beta<1$, %
then for any $t > 0$, $x \in \bR^d$ and $p \geq 2$,
\begin{equation}\label{E:UB-B12}
  \|u_1(t,x)\|_p \le c_{1} \exp\Big(c_{2}(p-1)^{\frac{1}{3-\alpha}}
  t^{\frac{4-\alpha-\beta}{3-\alpha}} \Big),
\end{equation}
where $c_1,c_2>0$ may depend on $d,\alpha,\beta$. For $d\ge3$ we use the
chaos-expansion definition of the solution in~\cite{balan.song:17:hyperbolic}.
In~\cite{balan:12:stochastic} the Riesz order is $d-\alpha$ in our notation. Its
Fourier estimate gives factorial decay $(n!)^{-(3-\alpha)/2}$ for the
$L^2(\Omega)$-norm of the $n$-th chaos component. Scaling in time and
hypercontractivity then bound the $L^p(\Omega)$-norm of the $n$-th chaos
component of $u_1(t,x)$ by
$\big(C(p-1)t^{4-\alpha-\beta}\big)^{n/2}/(n!)^{(3-\alpha)/2}$. Summing over $n$
using Stirling's formula gives~\eqref{E:UB-B12}. The hypercontractivity
inequality used here is
\begin{equation}\label{E:hyper}
  \|F\|_p \le (p-1)^{n/2} \|F\|_2 \quad \mbox{for any} \ F \in
  \cH_n
  \ \mbox{and}\ p\geq 2.
\end{equation}
Theorem 7.2 of~\cite{balan.conus:16:intermittency} gives the lower bound similar
to~\eqref{E:UB-B12} for $p = 2$ and $d \le 3$, based on a Feynman-Kac (FK)-type
representation for the second moment of the solution, which is inspired
by~\cite{dalang.mueller.ea:08:feynman-kac-type} and uses a planar Poisson
process (independent of the noise). This method works only in dimension $d \le
3$, and is an attempt to mimic the extraordinarily rich theory that exists for
the stochastic heat equation, based on the classical FK formula.

\medskip

In~\cite{CDST25}, the authors showed that for $d\in \{1,2,3\}$, the necessary and sufficient condition for the existence and uniqueness of the Skorohod solution of the hyperbolic Anderson model with Gaussian noise with covariance $\bE[\dot{W}(t,x)\dot{W}(s,y)]=|t-s|^{-\beta}\gamma(x-y)$ is
\begin{equation}
  \label{E:nec-suf1}
  \int_{\bR^d} \left( \frac{1}{1+|\xi|^2}\right)^{\frac{3-\beta}{2}} \mu(d\xi)<\infty,
\end{equation}
where $\beta \in (0,1)$ and $\gamma$ is Fourier transform of a tempered measure $\mu$, and proved that this condition is sufficient in the case $\beta=0$, i.e. when the noise is time-independent. If $\gamma$ satisfies the scaling property
\begin{equation}
  \label{E:scal-gam}
  \gamma(cx)=c^{-\alpha}\gamma(x) \quad \mbox{for any $c>0,x \in \bR^d$},
\end{equation}
for some $\alpha>0$, condition~\eqref{E:nec-suf1} becomes $\alpha+\beta<3$. The same problem for the Stratonovich solution in dimension $d\in\{1,2\}$ was treated in~\cite{CDST23}, where it was proved that the following condition is sufficient for the existence and uniqueness of the solution:
\begin{equation}
  \label{E:nec-suf2}
  \int_{\bR^d} \left( \frac{1}{1+|\xi|}\right)^{\rho_d-\beta-\eta} \mu(d\xi)<\infty, \quad \mbox{for some} \ \eta>0,
\end{equation}
where $\rho_1=1$ and $\rho_2=\frac{1}{2}$. If $\gamma$ satisfies~\eqref{E:scal-gam}, condition~\eqref{E:nec-suf2} becomes $\alpha+\beta<1$ if $d=1$ and $\alpha+\beta<\frac{1}{2}$ if $d=2$. Earlier, the first author proved in~\cite{balan22} the existence of the Stratonovich solution in the time-independent case, in dimension $d\in\{1,2\}$, under the condition:
\[
  \int_{\bR^d} \frac{1}{1+|\xi|} \mu(d\xi)<\infty.
\]

The exact asymptotic behaviour of the moments of the solution of the
hyperbolic Anderson model driven by a time-independent Gaussian noise was obtained by the first two authors and X. Chen in~\cite{balan.chen.ea:22:exact} in the Skorohod regime, and by X. Chen and Y.~Hu in~\cite{chen-hu} in the Stratonovich regime. These results are valid in any dimension $d\in \{1,2,3\}$.

In the time-dependent case, the picture is far from being complete.
In~\cite{chen:25:time-dependency}, X. Chen proved the large-time growth of the first moment of the Stratonovich solution, and proposed its sharp constant. But the problem of the exact asymptotic behaviour of the moments (paralleling the results of \cite{balan.chen.ea:22:exact,chen-hu} for the time-independent case mentioned above) has remained unsolved, for both Skorohod and Stratonovich regimes. The recent preprint~\cite{chen:26:stratonovich} by the second author gives a variational formula for the sharp constant governing the second moment, in the Stratonovich regime.

The goal of the present paper is to fill this gap and provide a complete solution to this problem in the Skorohod regime, in dimension $d=1$, for a Gaussian noise which has covariance function given by a product of two Riesz kernels, in space and time. Our matching upper and lower bounds will depend on a constant $D$ which is related to a variational problem for a deterministic wave equation, and governs both large-$t$ and large-$p$ asymptotics. As in \cite{balan.chen.ea:22:exact}, to treat simultaneously both asymptotics, we define
\begin{equation}\label{E:tp}
  t_{p} \coloneqq (p-1)^{\frac{1}{4-\alpha-\beta}} t, \qquad p
  \geq 2, \ t>0.
\end{equation}

Below is the main result of this article, which shows
that $\log\|u_\theta(t,x)\|_p\sim\log\|u_\theta(t_p,x)\|_2$ as $t_p\to\infty$.

\begin{theorem}\label{T:main-th}
  Let $t_p$ be given in~\eqref{E:tp}. There exists a constant
  $D \in (0,\infty)$ such that
  \begin{equation}\label{E:main-eq-tp}
    \lim_{t_p \to \infty} t_p^{-\frac{4-\alpha-\beta}{3-\alpha}} \log
    \|u_{\theta}(t,x)\|_p=\theta^{\frac{1}{3-\alpha}}D.
  \end{equation}
  Here and below, $t_p \to \infty$ means $t_p \to \infty$ jointly over $p\ge2$
  and $t > 0$. In particular,
  \begin{gather}\label{E:main-eq-p}
    \lim_{p \to \infty} p^{-\frac{4-\alpha}{3-\alpha}}
    \log \bE|u_{\theta}(t,x)|^p=\theta^{\frac{1}{3-\alpha}}
    t^{\frac{4-\alpha-\beta}{3-\alpha}} D, \quad \mbox{for any $t > 0$, and}\\
    \lim_{t \to \infty} t^{-\frac{4-\alpha-\beta}{3-\alpha}} \log
    \bE|u_{\theta}(t,x)|^p=p(p-1)^{\frac{1}{3-\alpha}}\theta^{\frac{1}{3-\alpha}}D,
    \quad \mbox{for any $p \geq 2$}.
    \label{E:main-eq-t}
  \end{gather}
\end{theorem}

\medskip

The constant $D$ in Theorem~\ref{T:main-th} is given below; its finiteness and
positivity are proved in Section~\ref{S:section-LB}:
\begin{equation}\label{E:def-D}
  D \coloneqq \sup_{f \in \cF_{+}}\left\{\Phi(V_f) -\frac{1}{2}\|f\|_{\cK}^2 \right\},
\end{equation}
with the following notation, which will be used throughout this article:

\medskip

\noindent\textit{(i)} $\cK \coloneqq L^2(E,\fm)$, where
$E \coloneqq \big((0,\infty) \times \bR\big)^2$ and
\[
  \fm (\ud a,\ud y,\ud b,\ud z) \coloneqq w_{\beta}(a) w_{\alpha}(b) \ud a\ud y
  \ud b\ud z \quad \mbox{with} \quad
  w_{\rho}(a) \coloneqq \frac{a^{\frac{\rho}{2}-1}}{\Gamma(\frac{\rho}{2}) }
  \ (a>0,\rho>0).
\]

\noindent\textit{(ii)} $V_f : [0,\infty) \times \bR \to [0,\infty)$ is {\em the
potential} of function  $f \in \cF_{+}\coloneqq \{f \in
C_c^{\infty}(E);f \geq 0\}$, defined by:
\[
  V_f(\xi) \coloneqq \int_{E} K_{\xi}(e)f(e) \fm(\ud e) \quad \mbox{for all} \ \xi=(t,x)
  \in [0,\infty) \times \bR,
\]
where, for all $t \geq 0$, $a, b > 0$ and $x, y, z, v \in \bR$,
\[
  K_{t,x}(a,y,b,z) \coloneqq \fg_a(t-y)\fg_b(x-z) \quad \mbox{with} \quad
  \fg_a(v)         \coloneqq \left(\frac{4a}{\pi}\right)^{\frac{1}{4}} e^{-2a v^2}.
\]
By Lemma~\ref{L:V-HypA} and the proof of Lemma~\ref{L:VN-HypA}, $V_f$ is bounded and smooth.

\medskip

\noindent\textit{(iii)} $\Phi$ is the {\em growth rate functional}, defined as follows: for
any bounded nonnegative measurable function $V: [0,1] \times \bR \to
[0,\infty)$,
\begin{equation}\label{E:def-Phi-V}
  \Phi(V) \coloneqq \int_0^1 \sqrt{\lambda_V(t)}\ud t,
\end{equation}
where
\begin{equation}\label{E:lambda-V}
  \lambda_V(t) \coloneqq \sup_{g \in H^1(\bR),\|g\|_{L^2(\bR)}=1}
  \left\{\int_{\bR}V(t,x)g^2(x)\ud x-\|g'\|_{L^2(\bR)}^2\right\},
\end{equation}
and $H^1(\bR)$ is the Sobolev space of order one (see~\cite[(1.4)--(1.5) and
Lemma~3.1]{balan.chen.ea:26:exact*1}). Note that the supremum defining $\lambda_V(t)$ can
be taken over a fixed countable $H^1$-dense subset of the $L^2$-unit sphere,
since its objective is continuous in the $H^1$-norm, if $V$ is bounded. Moreover, each term
is measurable in $t$, hence $\lambda_V$ is measurable. Under
Hypothesis~\ref{H:hypA} (given in Section~\ref{S:section-wave}), $\Phi(V)$ is
the exponential growth rate of the wave equation with potential $V$;
see~\cite[Theorem~1.1(i)]{balan.chen.ea:26:exact*1}.

\medskip

For the proof of the upper bound in
Theorem~\ref{T:main-th}, we use the moment comparison principle due to L\^e
(see~\cite{le:16:remark}, or Theorem B.1 of~\cite{balan.chen.ea:22:exact}),
which states that for any $p \geq 2$,
\begin{equation}\label{E:hyper2}
  \|u_{\theta}(t,x)\|_p \le \|u_{(p-1)\theta}(t,x)\|_2.
\end{equation}
This is a stronger form of hypercontractivity than \eqref{E:hyper}, which has
become standard in the literature, and reduces the upper bound problem to the
case $p = 2$.

The next step is to develop a link between the solution $u_{\theta}$
of~\eqref{E:HAM} and the solution $U_R^{Q}$ of a re-scaled wave equation
(equation
\eqref{E:HAM-R} below), which transforms the problem of super-exponential
growth of the moments of $u_{\theta}(t,x)$ with speed
$t_p^{\frac{4-\alpha-\beta}{3-\alpha}}$ to the exponential growth of the moments
of $U_R^Q(1,0)$; see Lemma~\ref{L:u-U} below. The upper index in $U_R^Q$
emphasizes the dependence on the covariance $Q$ of the noise:
\begin{equation}\label{E:def-Q}
  Q(t,x)=|t|^{-\beta}|x|^{-\alpha}, \quad t\ne0, x\in
  \bR\setminus\{0\}.
\end{equation}
(Formally, $\bE[\dot{W}(t,x)\dot{W}(s,y)] = Q(t-s,x-y)$. But $\dot{W}(t,x)$ does
not exist.) This trick eliminates the technicalities created by exponents, using
re-scaling. The core problem now remains to deal with $\bE|U_R^{Q}(1,0)|^2$, and
this is where a new idea emerges.

The idea starts with the observation that it is enough to
consider the analogue of equation \eqref{E:HAM-R} driven by a smooth Gaussian
noise $\dot{W}^{Q_N}$, whose covariance $Q_N$ is a truncation of $Q$. Unlike
$\dot{W}(t,x)$, the random field $\dot{W}^{Q_N}(t,x)$ exists, and it acts as a
potential for the wave equation.
The key property is that $\bE|U_R^{Q_N}(1,0)|^2$
coincides with the second moment of the solution of a deterministic wave
equation with {\em random} potential $Z_N/\sqrt{R}$, where $Z_N$ is a complex
random field whose real and imaginary parts are independent copies of
$\dot{W}^{Q_N}/\sqrt{2}$ (see Proposition~\ref{P:uQ-vR}). The core of the paper
consists in carefully importing into this stochastic framework an
exact asymptotic result which exists for the deterministic wave equation. (This
result is of independent interest, and is stated in Theorem~1.1(i)
of~\cite{balan.chen.ea:26:exact*1}; see Theorem~\ref{T:log-v-th} below.)
For a potential depending only on time and for a
comparison with the heat equation, see~\cite[Remarks~1.3 and~1.4,
respectively]{balan.chen.ea:26:exact*1}. For the importing procedure, we use the
Cameron-Martin space of the process $Z_N/\sqrt{R}$, when viewed as a random
element in the Banach space $\mathbf{B}$ of continuous functions on
$[0,1] \times [-R,R]$, and Borell-Sudakov-Cirel'son Theorem
from~\cite{ledoux:96:isoperimetry}.
The field $Z_N/\sqrt R$ has continuous sample paths by
Kolmogorov's criterion. The support of $G$ restricts the spatial variable to
$[-R,R]$ at $(1,0)$.

For the proof of the lower bound in Theorem~\ref{T:main-th},
we use a basic Malliavin calculus fact: if $F$ is a random variable with chaos
expansion $F = 1+\sum_{n \geq 1}I_n^{\fX}(h_n)$ with respect to an isonormal
Gaussian process $\{\fX(h)\}_{h \in \fH}$ and $Z_h =
\exp\big(\fX(h)-\frac{1}{2}\|h\|_{\fH}^2\big) =
\sum_{n \geq 0}\frac{1}{n!}I_n(h^{\otimes n})$, then by H\"older's
inequality,
\[
  |\bE[F Z_h]| \le \|F\|_p \|Z_h\|_q=\|F\|_p \,
  e^{\frac{q-1}{2}\|h\|_{\fH}^2}
  \quad \mbox{with} \ \frac{1}{p}+\frac{1}{q}=1.
\]
Therefore, by orthogonality of the Wiener chaos
spaces and the fact that $\fX(h) \sim
N(0,\|h\|_{\fH}^2)$,
\begin{equation}\label{E:LB-F}
  \|F\|_p \geq \Big(1+\sum_{n\geq 1}\langle \widetilde{h}_n,h^{\otimes
  n}\rangle_{\fH^{\otimes n}} \Big) e^{-\frac{1}{2(p-1)}
  \|h\|_{\fH}^2},
\end{equation}
where $\widetilde{h}_n$ is the symmetrization of $h_n$. This
reduces the lower bound problem to estimation of deterministic quantities
$\langle
\widetilde{h}_n,h^{\otimes n}\rangle_{\fH^{\otimes n}}$, which
are similar to those encountered for $\|F\|_2$.
In our case, we apply \eqref{E:LB-F} to $F =
U_R^{(p-1)^{-1}Q}(1,0)$, and $\fX = B$, the Gaussian white noise
on $(E,\fm)$, which is isometrically-linked to $W$, via Lemma~\ref{L:T-map}.
After some calculations, it turns out that $\bE[F Z_h]$ coincides with the
solution of the same deterministic wave equation with potential. The desired
lower bound follows from the limit in Theorem~1.1(i)
of~\cite{balan.chen.ea:26:exact*1}.

\medskip

The discussion above shows that at the core of the proof of Theorem~\ref{T:main-th} is a similar result about the exact asymptotic behaviour of the solution of a {\em deterministic} wave equation with potential, which the authors have proved in the companion paper~\cite{balan.chen.ea:26:exact*1}.

\medskip

This article is organized as follows. In Section~\ref{S:section-background}, we
present some preliminary material about Malliavin calculus and the solution. In
Section~\ref{S:section-rescaled}, we introduce the re-scaled equation and its
relationship with~\eqref{E:HAM}. Section~\ref{S:section-wave} presents the
connection between the equation with smooth noise and the deterministic wave
equation. Section~\ref{S:section-UB} proves the second-moment upper bound.
Sections~\ref{S:upper-p} and~\ref{S:section-LB} give the upper and lower
bounds, respectively, for general $p$.
Appendix \ref{S:appB} examines the Cameron-Martin
space. Appendix~\ref{S:static-comparison} gives the proof of
Proposition~\ref{P:static-comparison}.

\section{Background}\label{S:section-background}

In this section, we introduce the Malliavin
calculus preliminaries, and we present some background information about the
solution.

\medskip

Let $\cH$ be the completion of $\cD$ with respect to the inner product $\langle
\cdot,\cdot \rangle_{\cH}$. Then $W = \{W(\varphi);\varphi \in \cH \}$ is an
isonormal Gaussian process.

To emphasize the dependence on the covariance, we may denote $W = W^Q$.

We recall some background about Malliavin calculus. We refer the reader
to~\cite{nualart:06:malliavin} for more details. \medskip

\noindent\textbf{Wiener-chaos expansion} Let $\cF^W$ be the $\sigma$-field
generated by $W$. Any random variable $F \in L^2(\Omega,\cF^{W},\bP)$ has the
{\em Wiener-chaos expansion}:
\begin{equation}\label{E:chaos-F}
  F=\sum_{n\geq 0}I_n(f_n), \quad \mbox{for some} \ f_n \in \cH^{\otimes n},
\end{equation}
where $I_0(x) = x$, $f_0 = \bE(F)$, and for $n \geq 1$, $I_n :\cH^{\otimes n}
\to \cH_n$ is the multiple Wiener integral with respect to $W$, while $\cH_n$ is
the $n$-th Wiener chaos space of $W$, i.e. the closed linear span in
$L^2(\Omega)$ of $\{H_n(W(\varphi));\varphi \in \cH\}$. Here $H_n$ is the
Hermite polynomial:
\[
  H_n(x)=(-1)^n e^{\frac{x^2}{2}} \frac{\ud^n}{\ud x^n}(e^{-\frac{x^2}{2}})
  \quad \mbox{with} \quad
  e^{sx-\frac{s^2}{2}}=1+\sum_{n\geq1}\frac{s^n}{n!}H_n(x),
  \ s,x\in\bR.
\]
Then $I_n^W(\varphi^{\otimes n}) = H_n(W(\varphi))$ for any $\varphi \in \cH$
with $\|\varphi\|_{\cH} = 1$. Therefore,
\begin{equation}\label{E:In-Hn}
  I_n^W(\varphi^{\otimes n}) =\|\varphi\|_{\cH}^n H_n\left(
  \frac{W(\varphi)}{\|\varphi\|_{\cH}}\right) \quad \mbox{for any} \quad
  \varphi
  \in \cH.
\end{equation}
For any $F \in L^2(\Omega,\cF^{W},\bP)$, there is a unique expansion
\eqref{E:chaos-F} with symmetric kernels $f_n$.

The spaces $(\cH_n)_{n \geq 1}$ are orthogonal. For any $f \in \cH^{\otimes
n}$, $I_n(f) = I_n(\widetilde{f})$, where $\widetilde{f}$ is the symmetrization
of $f$, and
\[
  \bE|I_n(f)|^2 =n! \|\widetilde{f}\|_{\cH^{\otimes n}}^2.
\]

\noindent\textbf{Malliavin derivative} Let $\cS$ be the class of ``smooth''
random variables of the form $F = f(W(\varphi_1),\ldots,W(\varphi_n))$ for some
$f \in C_b^\infty(\bR^n)$ and $\varphi_1,\ldots,\varphi_n \in \cH$. The {\em
Malliavin derivative} of $F$ with respect to $W$ is defined as the following
element of $L^2(\Omega;\cH)$:
\[
  DF =\sum_{i=1}^{n}\frac{\partial}{\partial
  x_i}f(W(\varphi_1),\ldots,W(\varphi_n))\varphi_i.
\]
$D$ is closable, and its domain $\bD^{1,2}$ is the closure of $\cS$ with respect
to the norm:
\[
  \|F\|_{1,2}=\|F\|_2+\|DF\|_{L^2(\Omega;\cH)}.
\]
If $F \in \bD^{1,2}$ has the Wiener-chaos expansion \eqref{E:chaos-F} with
symmetric functions $f_n \in \cH^{\otimes n}$, then
\[
  D_{t,x}F=\sum_{n\geq 1}nI_{n-1}\big(f_n(\cdot,t,x) \big).
\]

\noindent\textbf{Skorohod integral} Let $\delta$ be the adjoint of $D$, whose
domain ${\rm Dom}(\delta)$ is the set of processes $u \in L^2(\Omega;\cH)$ for
which there exists $c > 0$ (depending on $u$), such that
\[
  |\bE\langle DF,u\rangle_{\cH} | \le c \|F\|_2 \quad \mbox{for all} \ F
  \in
  \bD^{1,2}.
\]
The following duality relation holds:
\[
  \bE\big(F \delta (u)\big)=\bE\langle DF,u\rangle_{\cH} \quad \mbox{for
  all}
  \ u \in {\rm Dom}(\delta) \ \mbox{and} \ F \in \bD^{1,2}.
\]
We say that $\delta$ is the {\em divergence operator} of $W$. For any $u \in
{\rm Dom}(\delta)$, we denote
\[
  \delta(u)=\int_{0}^{\infty} \int_{\bR}u(t,x)W(\delta t, \delta x),
\]
and we say that $\delta(u)$ is the {\em Skorohod integral} of $u$ with respect
to $W$.

\medskip

It is known from~\cite{balan:12:stochastic} that equation
\eqref{E:HAM}
has a unique solution $u_{\theta}$, and this solution has the Wiener-chaos
expansion:
\begin{equation}\label{E:chaos-u}
  u_{\theta}(t,x)=1+\sum_{n\geq 1} \theta^{n/2} I_n^W
  \big(f_{n}(\cdot,t,x)\big),
\end{equation}
where
\[
  f_n(t_1,x_1,\ldots,t_n,x_n,t,x)=\prod_{j=1}^{n}G(t_{j+1}-t_j,x_{j+1}-x_j)1_{\{0<t_1<\ldots<t_n<t\}},
\]
and $t_{n+1} = t,x_{n+1} = x$. Therefore,
\[
  \bE|u_{\theta}(t,x)|^2=1+\sum_{n\geq 1}\theta^n n!
  \|\widetilde{f}_{n}(\cdot,t,x)\|_{\cH^{\otimes n}}^2.
\]

By the spatial homogeneity of the noise $W$, $\{u_{\theta}(t,x)\}_{x \in
\bR}$
is strictly stationary. Therefore, it is enough to consider $x = 0$.

\section{Re-scaled equation}\label{S:section-rescaled}

In this section, we develop a connection between \eqref{E:HAM} and a re-scaled
equation, which will be used to prove \eqref{E:main-eq-tp}.

\medskip

More precisely, for any $R > 0$, we consider the equation:
\begin{align}
  \label{E:HAM-R}
  \begin{cases}
    \dfrac{\partial^2 U}{\partial t^2} (t,x)
    = R^2 \dfrac{\partial^2 U}{\partial x^2} (t,x)+R^{\frac{3}{2}} U(t,x)
    \dot{W}(t,x), \quad
    t>0, \ x \in \bR, \\[1em]
    U(0,x) = 1, \quad \dfrac{\partial U}{\partial t} (0,x)=0, \quad x \in \bR.
  \end{cases}
\end{align}

We say that a process $\{U_R^Q(t,x);(t,x) \in \bR_{+} \times \bR\}$ is a
solution of~\eqref{E:HAM-R} if it is $(\cF_t)$-adapted and satisfies:
\[
  U_{R}^Q(t,x)=1+R^{\frac{3}{2}}\int_0^t \int_{\bR} G_{R}(t-s,x-y)
  U_R^Q(s,y)
  W(\delta s,\delta y),
\]
where
\begin{equation}\label{E:def-GR}
  G_R(t,x)=\frac{1}{2R}1_{\{|x|<Rt\}}
\end{equation}
is the fundamental solution of $\frac{\partial^2 U}{\partial t^2} = R^2
\frac{\partial^2 U}{\partial x^2}$ on $\bR_{+} \times \bR$. Note that $G_1 =
G$. The superscript in $U_R^Q$ specifies the covariance
kernel~\eqref{E:def-Q} of the noise $\dot{W}$ in~\eqref{E:HAM-R}; recall that
$W = W^Q$. We also use~\eqref{E:HAM-R} with other covariance kernels, including
$(p-1)^{-1}Q$ in Lemma~\ref{L:u-U} and the truncation $Q_N$ in
Section~\ref{S:section-wave}.

Similarly to~\eqref{E:HAM}, equation \eqref{E:HAM-R} has a unique solution whose
chaos expansion is:
\begin{equation}\label{E:chaos-URQ}
  U_{R}^Q(t,x)=1+\sum_{n\geq 1} R^{\frac{3n}{2}} I_n^{W} \big(
  f_{n,R}(\cdot,t,x)\big),
\end{equation}
where
\begin{equation}\label{E:def-fnR}
  f_{n,R}(t_1,x_1,\ldots,t_n,x_n,t,x)=\prod_{j=1}^{n}G_R(t_{j+1}-t_j,x_{j+1}-x_j)1_{\{0<t_1<\ldots<t_n<t\}}.
\end{equation}
Therefore,
\begin{equation}\label{E:2m-URQ}
  \bE|U_{R}^Q(t,x)|^2=1+\sum_{n\geq 1} R^{3n} n!
  \|\widetilde{f}_{n,R}(\cdot,t,x)\|_{\cH^{\otimes n}}^2.
\end{equation}

\begin{lemma}
  \label{L:u-U}
  For any $t > 0$, $x \in \bR$ and $p \geq 2$, with value $t_p$
  given in~\eqref{E:tp},
  \[
    u_{\theta}(t,x) \stackrel{d}{=} U_{R_{p,t}}^{(p-1)^{-1}Q}(1,0) \quad \mbox{with} \quad
    R_{p,t} \coloneqq \theta^{\frac{1}{3-\alpha}} t_p^{\frac{4-\alpha-\beta}{3-\alpha}}
    = \big((p-1) \theta\big)^{\frac{1}{3-\alpha}} t^{\frac{4-\alpha-\beta}{3-\alpha}},
  \]
  where $\stackrel{d}{=}$ denotes equality in distribution and
  $U_{R}^{(p-1)^{-1}Q}$ is the solution of~\eqref{E:HAM-R} with $\dot{W}$
  replaced by $(p-1)^{-1/2}\dot{W}$.
  Consequently,
  with the same $p$ on both sides along each sequence and $R=R_{p,t}$, we
  have:
  \begin{equation}\label{E:u-U-sup}
    \limsup_{t_p \to \infty} t_p^{-\frac{4-\alpha-\beta}{3-\alpha}} \log \|
    u_{\theta}(t,x)\|_p =\theta^{\frac{1}{3-\alpha}} \limsup_{R \to \infty}
    \frac{1}{R}\log \|U_{R}^{(p-1)^{-1}Q}(1,0)\|_p,%
  \end{equation}
  and
  \begin{equation}\label{E:u-U-inf}
    \liminf_{t_p \to \infty} t_p^{-\frac{4-\alpha-\beta}{3-\alpha}} \log \|
    u_{\theta}(t,x)\|_p =\theta^{\frac{1}{3-\alpha}} \liminf_{R \to \infty}
    \frac{1}{R}\log \|U_{R}^{(p-1)^{-1}Q}(1,0)\|_p.%
  \end{equation}
\end{lemma}

\begin{proof}
  We fix $(t,x) \in (0,\infty) \times \bR$ and $R > 0$. For any $\varphi \in
  \cD$, we define
  \begin{equation}\label{E:def-J}
    (J \varphi)(r,z)=t^{-\frac{4-\alpha-\beta}{2}} R^{1-\frac{\alpha}{2}}
    \varphi \left( \frac{r}{t},\frac{R(z-x)}{t}\right), \quad \mbox{for}
    \quad
    r>0,z\in \bR.
  \end{equation}
  By scaling, $\langle J \varphi, J\psi \rangle_{\cH} = \langle \varphi,
  \psi
  \rangle_{\cH}$. Therefore, if we define
  \[
    \widehat{W} (\varphi)=W(J \varphi) \quad \mbox{for all} \quad \varphi
    \in
    \cD,
  \]
  then the map $\cD \ni \varphi \mapsto \widehat{W}(\varphi) \in L^2(\Omega)$ is
  an isometry which can be extended to $\cH$. The processes $W$ and
  $\widehat{W}$ have the same law.

  Note that $G_R$ has the following scaling property: for any $t > 0$ and $x
  \in \bR$,
  \[
    G_R(t,x)=\frac{1}{R} G\left(at,\frac{ax}{R}\right) \quad \mbox{for any}
    \ a>0.
  \]
  Therefore, for any $(r_1,z_1),\ldots,(r_n,z_n) \in (0,\infty) \times \bR$,
  \begin{equation}\label{E:fnR-scale}
    f_{n,R} \left(\frac{r_1}{t},\frac{R(z_1-x)}{t}, \ldots,
    \frac{r_n}{t},\frac{R(z_n-x)}{t},1,0\right)=\frac{1}{R^n}f_n(r_1,z_1,\ldots,r_n,z_n,t,x),
  \end{equation}
  and consequently, by approximation,
  \[
    J^{\otimes n} f_{n,R}(\cdot,1,0)=t^{-\frac{4-\alpha-\beta}{2}n}
    R^{-\frac{\alpha n}{2}}f_n(\cdot,t,x).
  \]

  Indeed, the bounded, compactly supported kernel has discontinuities only
  on finitely many hyperplanes. Mollification and cutoff give uniformly bounded
  smooth approximants with common compact support converging almost everywhere.
  The product covariance is locally integrable, so dominated convergence gives
  convergence in the tensor norm before and after the change of variables.
  Smooth compactly supported functions may in turn be approximated uniformly
  on a common compact rectangle by finite sums of smooth product functions.

  Moreover, using~\eqref{E:In-Hn}, we have:
  \begin{equation}\label{E:IWW-hat}
    I_n^{W}(J^{\otimes n} h)=I_n^{\widehat{W}}(h) \quad \mbox{for all} \quad
    h
    \in \cH^{\otimes n}.
  \end{equation}

  Indeed, this relation holds for $h = \varphi ^{\otimes n}$ and $\varphi
  \in \cH$, since
  \[
    I_n^{W}(J^{\otimes n} \varphi ^{\otimes n})=I_n^{W}\big( (J \varphi
    )^{\otimes n} \big)=\| J \varphi \|_{\cH}^n \, H_n \left( \frac{W(J
    \varphi
    )}{\|J\varphi \|_{\cH}}\right)=\|\varphi \|_{\cH}^n H_n\left(
    \frac{\widehat{W}(\varphi )}{\|\varphi
    \|_{\cH}}\right)=I_n^{\widehat{W}}(\varphi ^{\otimes n}),
  \]
  and can be extended to any $h \in \cH^{\otimes n}$. Therefore,
  using~\eqref{E:chaos-u}, we have:
  \begin{align}
    \nonumber
    u_{\theta}(t,x)& =1+\sum_{n\geq 1} \theta^{\frac{n}{2}}
    t^{\frac{4-\alpha-\beta}{2}n} R^{\frac{\alpha
    n}{2}}I_n^{W}\big(J^{\otimes
    n}f_{n,R}(\cdot,1,0)\big) \\
    \label{E:u-rep1}
    &=1+\sum_{n\geq 1} \big( \theta t^{4-\alpha-\beta} R^{\alpha}
    \big)^{\frac{n}{2}} I_n^{\widehat{W}}\big(f_{n,R}(\cdot,1,0)\big)\\
    \nonumber
    &=1+\sum_{n\geq 1} \big( (p-1)\theta t^{4-\alpha-\beta} R^{\alpha}
    \big)^{\frac{n}{2}}
    (p-1)^{-\frac{n}{2}}I_n^{\widehat{W}}\big(f_{n,R}(\cdot,1,0)\big),
  \end{align}
  for any $R > 0$. Since $R_{p,t}^{3-\alpha} = (p-1)\theta
  t^{4-\alpha-\beta}$, taking $R = R_{p,t}$, we get:
  \[
    u_{\theta}(t,x)=1+\sum_{n\geq 1} R_{p,t}^{\frac{3n}{2}}
    (p-1)^{-\frac{n}{2}}I_n^{\widehat{W}}\big(f_{n,R}(\cdot,1,0)\big).
  \]

  On the other hand, by~\eqref{E:chaos-URQ}, $U_{R}^{(p-1)^{-1}Q}(1,0)$ has the
  chaos expansion:
  \begin{equation}\label{E:chaos-URQp}
    U_{R}^{(p-1)^{-1}Q}(1,0)=1+\sum_{n\geq 1} R^{\frac{3n}{2}}
    (p-1)^{-\frac{n}{2}}I_n^{W}\big(f_{n,R}(\cdot,1,0)\big),
  \end{equation}
  for any $R > 0$. When $R = R_{p,t}$, the last series has the same distribution
  as $u_{\theta}(t,x)$, since $W$ and $\widehat{W}$ have the same distribution.

  Finally, for fixed $\theta>0$, $R_{p,t}\to\infty$ if and only if
  $t_p\to\infty$, and $t_p^{-\frac{4-\alpha-\beta}{3-\alpha}} =
  \theta^{\frac{1}{3-\alpha}}/R_{p,t}$. Taking the limsup and the liminf of the
  equal log norms along corresponding sequences with the same $p$
  gives~\eqref{E:u-U-sup} and~\eqref{E:u-U-inf}.
  This proves Lemma~\ref{L:u-U}.
\end{proof}

\section{Connection with a deterministic wave equation}\label{S:section-wave}

In this section, we consider equation~\eqref{E:HAM-R} with a smooth noise
$\dot{W}^{Q_N}$ replacing $\dot{W}$, and we show that this equation is connected
with a deterministic wave equation.

\medskip

\subsection{Motivation}\label{SS:motivation}

In this section, we give the motivation behind the definition of smooth noise
$\dot{W}^{Q_N}$.

Note that for any $\rho > 0$,
\[
  h_{\rho}(v)\coloneqq  \int_{0}^{\infty}e^{-av^2} w_{\rho}(a)\ud a =|v|^{-\rho} \quad
  \mbox{for all} \ v \in \bR\setminus\{0\}.
\]
Therefore, $Q(t,x) = h_{\beta}(t)h_{\alpha}(x)$ for $t\ne0$ and
$x \in \bR\setminus\{0\}$.

Moreover, from the semigroup property of the heat kernel, it follows that for
any $a > 0$
\begin{equation}\label{E:semigroup}
  (\fg_{a}*\fg_a)(v)=e^{-a v^2} \quad \mbox{for any} \ v \in \bR,
\end{equation}
and therefore, for any $\xi = (t,x),\xi' = (t',x') \in [0,\infty) \times
\bR$
with $t\not = t'$ and $x\not = x'$,
\begin{equation}\label{E:KK-int}
  \int_{E}K_{\xi}(e)K_{\xi'}(e) \fm(
  \ud e)=h_{\beta}(t-t')h_{\alpha}(x-x')=Q(t-t',x-x').
\end{equation}
But $K_{\xi} \not \in \cK$, since $\|K_{\xi}\|_{\cK}^2 = Q(0,0) = \infty$.

\begin{lemma}\label{L:T-map}
  For any $\varphi \in \cD$, we define the element $T \varphi \in \cK$ by:
  \[
    (T \varphi)(e)=\int_0^{\infty}\int_{\bR} K_{t,x}(e)\varphi(t,x) \ud t
    \ud x,
    \quad \mbox{for all} \ e \in E.
  \]
  Then $T$ is an isometry between $\cD$ and $\cK$ which can be extended to
  $\cH$, and we can take $$W(\varphi)=B(T\varphi) \quad \mbox{for all}
  \ \varphi
  \in \cH,$$ where $B$ is a Gaussian white noise on $(E, \fm)$, i.e. $B =
  \{B(\phi);\phi \in \cK\}$ is a zero-mean Gaussian process with covariance
  $\bE[B(\phi)B(\psi)] = \langle \phi,\psi \rangle_{\cK}$ for all $\phi,\psi
  \in \cK$.
\end{lemma}

\begin{proof}
  Using Fubini theorem and~\eqref{E:KK-int}, we have: for any $\varphi,\psi \in
  \cD$,
  \begin{align*}
    \langle T \varphi, T \psi \rangle_{\cK}&=\int_{E}
    (T\varphi)(e)(T\psi)(e)
    \fm(\ud e)\\
    &=\int_{([0,\infty)\times \bR)^2} \left(\int_{E} K_{\xi}(e) K_{\xi'}(e)
    \fm(\ud e)\right) \varphi(\xi) \psi(\xi') \ud \xi \ud \xi'\\
    &=\int_{([0,\infty)\times \bR)^2} Q(\xi-\xi') \varphi(\xi)
    \psi(\xi') \ud \xi \ud \xi'.
  \end{align*}
  On the closed image $T(\cH)$ define $B(T\varphi)=W(\varphi)$. On a product
  extension take an independent isonormal process on its orthogonal complement
  in $\cK$. Their sum defines $B$ on all of $\cK$ without changing the law of
  the original noise or its solution.
  This proves Lemma~\ref{L:T-map}.
\end{proof}

Throughout, an integral in space-time variables written without a domain is
taken over $[0,\infty)\times\bR$ in each variable.

\subsection{Equation with smooth noise}

Motivated by the discussion above, we now define, for any integer
$N \geq 2$,
\[
  Q_N(t,x)=h_{\beta,N}(t)h_{\alpha,N}(x) \quad \mbox{for all $t,x\in\bR$},
\]
where
\[
  h_{\rho,N}(v)\coloneqq  \int_{\frac{1}{N}}^{N}e^{-av^2} w_{\rho}(a)\ud a \quad
  \mbox{for all} \ v \in \bR.
\]

Let $\cK_N = L^2(E_N,\fm)$, where $E_N = \Big([\frac{1}{N},N] \times
\bR\Big)^2$. Similarly to~\eqref{E:KK-int}, using relation
\eqref{E:semigroup}, we
obtain that for any $\xi = (t,x),\xi' = (t',x') \in [0,\infty) \times \bR$,
\[
  \int_{E_N}K_{\xi}(e)K_{\xi'}(e) \fm(
  \ud e)=h_{\beta,N}(t-t')h_{\alpha,N}(x-x')=Q_N(t-t',x-x').
\]
But $K_{\xi} \in \cK_N$, since $\|K_{\xi}\|_{\cK_N}^2 = Q_N(0,0) < \infty$.
Moreover, $Q_N$ is non-negative definite: for any $n \geq 1$ and for any
$\xi_1,\ldots,\xi_n \in [0,\infty) \times \bR$ and $a_1,\ldots,a_n \in \bR$,
\[
  \sum_{i,j=1}^n a_i a_j Q_N(\xi_i-\xi_j)=\sum_{i,j=1}^n a_i a_j \langle
  K_{\xi_i},K_{\xi_j}\rangle_{\cK_N}=\| \sum_{i=1}^n a_i
  K_{\xi_i}\|_{\cK_N}^2
  \geq 0.
\]

Hence, there exists a zero-mean Gaussian random field
$\{\dot{W}^{Q_N}(t,x);t \geq 0,x \in \bR\}$ with
\[
  \bE[\dot{W}^{Q_N}(t,x)\dot{W}^{Q_N}(s,y)]=Q_{N}(t-s,x-y)
\]
for all $t \geq 0$, $s \geq 0$ and $x,y \in \bR$. We take
\begin{equation}\label{E:def-WQN}
  \dot{W}^{Q_N}(t,x)=B_N(K_{t,x}) \quad \mbox{for all} \ t\geq 0,x\in \bR,
\end{equation}
where $B_N$ is a Gaussian white noise on $(E_N, \fm)$, i.e. $B_N =
\{B_N(\phi);\phi \in \cK_N\}$ is a zero-mean Gaussian process with
covariance $\bE[B_N(\phi)B_N(\psi)] = \langle \phi,\psi \rangle_{\cK_N}$ for all
$\phi,\psi \in \cK_N$.

\medskip

We define
\[
  W^{Q_N}(\varphi)=\int_0^{\infty}\int_{\bR}\varphi(t,x)\dot{W}^{Q_N}(t,x)
  \ud t\ud x, \quad \mbox{for all} \ \varphi \in \cD.
\]
Then for any $\varphi,\psi \in \cD$,
\[
  \bE[W^{Q_N}(\varphi)W^{Q_N}(\psi)]=\int_{((0,\infty) \times \bR)^2}
  \varphi(t,x)\psi(s,y)Q_N(t-s,x-y)\ud t\ud x \ud s\ud y=:\langle \varphi,\psi
  \rangle_{\cH^{Q_N}}.
\]
Let $\cH^{Q_N}$ be the completion of $\cD$ with respect to $\langle
\cdot,\cdot
\rangle_{\cH^{Q_N}}$. The map $\varphi \mapsto W^{Q_N}(\varphi) \in
L^2(\Omega)$ is an isometry which can be extended to $\cH^{Q_N}$. Then $W^{Q_N}
=
\{W^{Q_N}(\varphi);\varphi \in \cH^{Q_N}\}$ is an isonormal Gaussian
process.

\medskip

Let $U_R^{Q_N}$ be the solution of equation \eqref{E:HAM-R} with $\dot{W}$
replaced by $\dot{W}^{Q_N}$, i.e. $\{U_R^{Q_N}(t,x); t \geq 0,x \in \bR\}$ is
adapted to the filtration generated by $W^{Q_N}$ and satisfies
\[
  U_{R}^{Q_N}(t,x)=1+R^{\frac{3}{2}}\int_0^t \int_{\bR} G_{R}(t-s,x-y)
  U_R^{Q_N}(s,y) W^{Q_N}(\delta s, \delta y).
\]
Similarly to~\eqref{E:chaos-URQ}, this solution has the chaos expansion:
\[
  U_{R}^{Q_N}(t,x)=1+\sum_{n\geq 1} R^{\frac{3n}{2}}
  I_{n}^{W^{Q_N}}(f_{n,R}(\cdot,t,x)),
\]
where %
$I_n^{W^{Q_N}}$ is the multiple Wiener integral of order $n$ with respect to
$W^{Q_N}$.

\begin{lemma}\label{L:mom2-uQ}
  For any $t \geq 0$ and $x \in \bR$,
  \[
    \bE|U_{R}^{Q_N}(t,x)|^2=1+\sum_{n\geq 1} R^{3n} \int f_{n,R}(\pmb{\xi},t,x)f_{n,R}(\pmb{\eta},t,x)
    \sum_{\sigma \in S_n} \prod_{j=1}^{n} Q_N(\xi_j-\eta_{\sigma(j)})
    \ud \pmb{\xi} \ud \pmb{\eta},
  \]
  where $\pmb{\xi} = (\xi_1,\ldots,\xi_n)$, $\pmb{\eta} =
  (\eta_1,\ldots,\eta_n)$, $S_n$ is the set of permutations of $1,\ldots,n$,
  and, by the convention at the end of Section~\ref{SS:motivation}, each
  $\xi_j$ and $\eta_j$ is integrated over $[0,\infty)\times\bR$.
\end{lemma}

\begin{proof}
  Similarly to~\eqref{E:2m-URQ},
  \[
    \bE|U_{R}^{Q_N}(t,x)|^2=1+\sum_{n\geq 1} R^{3n}
    n!\|\widetilde{f}_{n,R}(\cdot,t,x)\|_{(\cH^{Q_N})^{\otimes n}}^2.
  \]
  We will use the fact that
  $\|\widetilde{f}\|_{(\cH^{Q_N})^{\otimes n}}^2 = \langle
  f,\widetilde{f}\rangle_{(\cH^{Q_N})^{\otimes n}}$ for any $f \in
  (\cH^{Q_N})^{\otimes n}$. Hence,
  \begin{align*}
    & \|\widetilde{f}_{n,R}(\cdot,t,x)\|_{(\cH^{Q_N})^{\otimes n}}^2
    =\langle
    f_{n,R}(\cdot,t,x),\widetilde{f_{n,R}}(\cdot,t,x)\rangle_{(\cH^{Q_N})^{\otimes
    n}}\\
    & =\frac{1}{n!} \sum_{\sigma \in S_n} \int f_{n,R}(\xi_1,\ldots,\xi_n,t,x)
    f_{n,R}(u_{\sigma(1)},\ldots,u_{\sigma(n)},t,x)
    \prod_{j=1}^{n}Q_N(\xi_j-u_j) \ud \pmb{\xi} \ud \pmb{u} \\
    & =\frac{1}{n!} \sum_{\sigma \in S_n} \int f_{n,R}(\xi_1,\ldots,\xi_n,t,x)
    f_{n,R}(\eta_{1},\ldots,\eta_{n},t,x)
    \prod_{j=1}^{n}Q_N(\xi_j-\eta_{\sigma^{-1}(j)}) \ud \pmb{\xi}
    \ud \pmb{\eta}.
  \end{align*}
  The conclusion of Lemma~\ref{L:mom2-uQ} follows.
\end{proof}

\medskip

Let $X_N$ and $Y_N$ be two independent copies of $\dot{W}^{Q_N}$, and define
\begin{equation}\label{E:def-ZN}
  Z_N=\frac{X_N+iY_N}{\sqrt{2}}.
\end{equation}
Note that for any $\xi,\eta \in [0,\infty) \times \bR$,
\begin{equation}\label{E:EZZ}
  \bE[Z_N(\xi)Z_N(\eta)]=0 \quad \mbox{and} \quad \bE[Z_N(\xi)
  \overline{Z_N(\eta)}]=Q_N(\xi-\eta).
\end{equation}

We consider the $\bC$-valued random
field $Z_N = \{Z_N(t,x);t \in [0,1],x \in \bR\}$, and we let
$v_R\big(\frac{Z_N}{\sqrt{R}};\cdot\big)$ be the (pathwise) solution of the
deterministic wave equation \eqref{E:HAM-v} with potential $V = Z_N/\sqrt{R}$.
At each point we define this solution by the series~\eqref{E:chaos-v} on its
backward cone, where the following bound gives absolute convergence.
For $\xi,\eta$ in a compact rectangle, \eqref{E:def-Cjk} and
Lemma~\ref{L:deriv-XN} bound $\|Z_N(\xi)-Z_N(\eta)\|_2$ by a constant times
$|\xi-\eta|$. For Gaussian increments the same holds in $L^p(\Omega)$, $p > 2$.
Kolmogorov's criterion gives continuous versions on rectangles, agreeing on
overlaps almost surely by continuity on a countable dense set. For this version,
on the cone ending at $(t,x)$ let $M = \sup|Z_N|/\sqrt R$.
Then $M < \infty$ and, by
Remark~\ref{R:support-kernel}, the mass calculation in the proof
of~\cite[Lemma~2.1]{balan.chen.ea:26:exact*1} bounds the absolute sum
in~\eqref{E:chaos-v} by $\sum_{n\ge0}(Rt\sqrt M)^{2n}/(2n)!$. Then
\[
  v_R\Big(\frac{Z_N}{\sqrt{R}};t,x\Big)=1+\sum_{n\geq 1} R^{\frac{3n}{2}}
  J_{n,R}(t,x) \quad \mbox{for} \ (t,x) \in [0,1] \times \bR,
\]
where
\[
  J_{n,R}(t,x)=\int_{((0,\infty) \times \bR)^n}
  f_{n,R}(\xi_1,\ldots,\xi_n,t,x)
  \prod_{j=1}^{n} Z_N(\xi_j) \ud \xi_1 \ldots \ud \xi_n.
\]

For $n\ge1$, multiplying $Z_N$ by
$e^{i\pi/n}$ preserves its joint law
by~\eqref{E:def-ZN} and multiplies $J_{n,R}(t,x)$ by $-1$.
Its expectation exists by H\"older's inequality, the finite Gaussian
moments and the finite mass of $f_{n,R}$. Hence, $\bE[J_{n,R}(t,x)] = 0$.

\begin{lemma}\label{L:SecondMom}
  For any $\xi = (t,x) \in [0,\infty) \times
  \bR$,
  and $R > 0$,
  \begin{equation}\label{E:Jn-orthog}
    \bE[J_{n,R}(\xi) \overline{J_{m,R}(\xi) }]=0 \quad \mbox{if} \quad
    n\not=m,
  \end{equation}
  and
  \begin{equation}\label{E:2m-JnR}
    \bE|J_{n,R}(\xi)|^2 =\int
    f_{n,R}(\pmb{\xi},t,x)f_{n,R}(\pmb{\eta},t,x) \sum_{\sigma \in S_n}
    \prod_{j=1}^{n} Q_N(\xi_j-\eta_{\sigma(j)}) \ud \pmb{\xi}
    \ud \pmb{\eta}.
  \end{equation}
\end{lemma}

\begin{proof}
  The same H\"older bound on the Gaussian products and the finite mass of both
  kernels justify Fubini's theorem. Thus,
  \[
    \begin{aligned}
      \bE[J_{n,R}(\xi) \overline{J_{m,R}(\xi) }]&=\int_{((0,\infty) \times
      \bR)^{n}} \int_{((0,\infty) \times \bR)^{m}}f_{n,R}(\pmb{\xi},\xi)
      f_{m,R}(\pmb{\eta},\xi)\\ &\qquad{}\times\bE\left[\prod_{j=1}^{n}
      Z_N(\xi_j) \prod_{k=1}^{m} \overline{Z_N(\eta_k)}\right] \ud \pmb{\xi}
      \ud \pmb{\eta},
    \end{aligned}
  \]
  where $\pmb{\xi} = (\xi_1,\ldots,\xi_n)$ and $\pmb{\eta} =
  (\eta_1,\ldots,\eta_m)$. First, we will prove that if $n\not = m$,
  \[
    \bE\left[\prod_{j=1}^{n} Z_N(\xi_j) \prod_{k=1}^{m}
    \overline{Z_N(\eta_k)}\right] =0 \quad \mbox{for all} \quad
    (\pmb{\xi},\pmb{\eta})\in ((0,\infty)\times \bR)^{n+m}.
  \]

  We fix such a pair $(\pmb{\xi},\pmb{\eta})$. We consider the moment
  generating function of the Gaussian vector $\big(Z_N(\xi_1),
  \ldots,Z_N(\xi_n), \overline{Z_N(\eta_1)}, \ldots,
  \overline{Z_N(\eta_m)}\big)$:
  \[
    M(\pmb{t},\pmb{s})=\bE\Big[ \exp\Big(\sum_{j=1}^n t_j Z_N(\xi_j) +
    \sum_{k=1}^m s_k\overline{Z_N(\eta_k)} \Big)\Big]
  \]
  for any $\pmb{t} = (t_1,\ldots,t_n) \in \bR^n$ and $\pmb{s} =
  (s_1,\ldots,s_m) \in \bR^m$. On one hand,
  \[
    \frac{\partial^{n+m}M}{\partial t_1 \ldots \partial t_n \partial s_1
    \ldots
    \partial s_m}(\pmb{0},\pmb{0})=
    \bE\left[\prod_{j=1}^{n} Z_N(\xi_j) \prod_{k=1}^{m}
    \overline{Z_N(\eta_k)}\right].
  \]
  On the other hand, by Gaussianity, $M(\pmb{t},\pmb{s}) = \exp\big(\sum_{j =
  1}^n \sum_{k = 1}^m t_j s_k Q_N(\xi_j-\eta_k)\big)$, and hence (because
  $n\not = m$),
  \[
    \frac{\partial^{n+m}M}{\partial t_1 \ldots \partial t_n \partial s_1
    \ldots
    \partial s_m}(\pmb{0},\pmb{0})=0.
  \]

  For the second statement, we will use Wick's formula for the zero-mean
  Gaussian vector $(G_1,\ldots,G_{2n}) =
  \big(Z_N(\xi_1),\ldots,Z_N(\xi_n),\overline{Z_N(\eta_1)},\ldots,\overline{Z_N(\eta_n)}\big)$:
  \[
    \bE[G_1 \ldots G_{2n}] =\sum_{\pi} \prod_{\ell=1}^{n} \bE[G_{i_{\ell}}
    G_{j_{\ell}}],
  \]
  where the sum is over all partitions $\pi =
  \{\{i_1,j_1\},\ldots,\{i_n,j_n\}\}$ of $\{1,\ldots,2n\}$ consisting of blocks
  of 2 elements.

  By~\eqref{E:EZZ}, if $i_{\ell},j_{\ell} \in \{1,\ldots,n\}$ or
  $i_{\ell},j_{\ell} \in \{n+1,\ldots,2n\}$, then $\bE[G_{i_{\ell}}G_{j_{\ell}}]
  = 0$, whereas if $i_{\ell} \in \{1,\ldots,n\}$ and $j_{\ell} \in
  \{n+1,\ldots,2n\}$, then
  \[
    \bE[G_{i_{\ell}}G_{j_{\ell}}]=Q_N(\xi_j-\eta_k),
  \]
  where $G_{i_{\ell}} = Z_N(\xi_j)$ and $G_{j_{\ell}} = \overline{Z_N(\eta_k)}$.
  Note that there is a bijection between the set of partitions $\pi =
  \{\{i_1,j_1\},\ldots,\{i_n,j_n\}\}$ with $i_{\ell} \in \{1,\ldots,n\}$ and
  $j_{\ell} \in \{n+1,\ldots,2n\}$, and the set $\bigl\{\{\{j,n+\sigma(j)\}:1
  \le j \le n\}:\sigma \in S_n\bigr\}$. Hence,
  \[
    \bE\left[\prod_{j=1}^{n} Z_N(\xi_j) \prod_{k=1}^{n}
    \overline{Z_N(\eta_k)}\right]=\sum_{\sigma \in S_n} \prod_{j=1}^{n}
    Q_N(\xi_j-\eta_{\sigma(j)}).
  \]
  The conclusion of Lemma~\ref{L:SecondMom} follows.
\end{proof}

\begin{proposition}\label{P:uQ-vR}
  For any $(t,x) \in [0,1] \times \bR$, we have:
  \[
    \bE|U_R^{Q_N}(t,x)|^2=\bE\Big|v_R\Big(\frac{Z_N}{\sqrt{R}};t,x\Big)\Big|^2.
  \]
\end{proposition}

\begin{proof} By~\eqref{E:2m-JnR}, $|Q_N| \le
  Q_N(0,0)$ and the mass of $f_{n,R}$,
  \[
    \bE|J_{n,R}(t,x)|^2\le n!\,Q_N(0,0)^n
    \left(\frac{t^{2n}}{(2n)!}\right)^2.
  \]
  This is summable after multiplication by $R^{3n}$. By orthogonality
  the series converges in $L^2(\Omega)$, with the same limit as its
  pathwise sum. Using~\eqref{E:Jn-orthog} and the fact that
  $\bE[J_{n,R}(\xi)] = 0$, we have:
  \[
    \bE\Big|v_R\Big(\frac{Z_N}{\sqrt{R}};t,x\Big)\Big|^2=1+\sum_{n\geq
    1}R^{3n}
    \bE|J_{n,R}(t,x)|^2.
  \]
  The conclusion follows by~\eqref{E:2m-JnR} and Lemma~\ref{L:mom2-uQ}.
  This proves Proposition~\ref{P:uQ-vR}.
\end{proof}

Let $V:[0,1]\times\bR \to \bC$ be a measurable function. We consider the
deterministic wave equation with potential $V$
(see~\cite[(1.1)--(1.2)]{balan.chen.ea:26:exact*1}):
\begin{align}
  \label{E:HAM-v}
  \begin{cases}
    \dfrac{\partial^2 v}{\partial t^2}(t,x)
    = R^2\dfrac{\partial^2 v}{\partial x^2}(t,x)
    +R^2v(t,x)V(t,x),
    \quad t\in[0,1],\ x\in\bR,\\[1em]
    v(0,x)=1,\quad \dfrac{\partial v}{\partial t}(0,x)=0,
    \quad x\in\bR.
  \end{cases}
\end{align}

By Duhamel's principle, a solution $v$ of~\eqref{E:HAM-v} satisfies the integral
equation:
\begin{equation}\label{E:duhamel-v}
  \begin{aligned}
    v(t,x)
    &=1+R^2\int_0^t\int_{\bR}
    G_R(t-s,x-y)v(s,y)V(s,y)\,\ud y\,\ud s\\
    &=1+\frac R2\int_0^t\int_{x-R(t-s)}^{x+R(t-s)}
    V(s,y)v(s,y)\,\ud y\,\ud s,
  \end{aligned}
\end{equation}
where $G_R$ is given by~\eqref{E:def-GR}.

By~\cite[Lemma~2.1, (2.2)]{balan.chen.ea:26:exact*1}, if $V$ is bounded,
then~\eqref{E:HAM-v} has a unique bounded mild solution $v_R(V;\cdot)$. For
every $t \in [0,1]$ and $x \in \bR$, this solution is given by the absolutely
and uniformly convergent series
\begin{equation}\label{E:chaos-v}
  v_R(V;t,x)
  =1+\sum_{n\geq1}R^{2n}
  \int_{([0,1]\times\bR)^n}
  f_{n,R}(\xi_1,\ldots,\xi_n,t,x)
  \prod_{j=1}^n V(\xi_j)\,\ud \xi_1\ldots \ud \xi_n,
\end{equation}
where $f_{n,R}(\cdot,t,x)$ is given by~\eqref{E:def-fnR}.

\begin{remark}\label{R:support-kernel}
  {\rm By~\cite[Remark~2.1, (2.4)]{balan.chen.ea:26:exact*1}, on the support of
  $f_{n,R}$ each $\xi_j = (s_j,x_j)$ satisfies $|x-x_j| \le R(t-s_j)$.
  Consequently, $f_{n,R}(\cdot,t,x)$ is supported in
  $\bigl([0,t]\times[x-Rt,x+Rt]\bigr)^n$. In particular, $f_{n,R}(\cdot,1,0)$ is
  supported in $\bigl([0,1]\times[-R,R]\bigr)^n$. Consequently,
  \begin{equation}\label{E:int-CR}
    v_R(V;1,0)=1+\sum_{n\ge 1} R^{2n}\int_{([0,1]\times[-R,R])^n}
    f_{n,R}(\xi_1,\ldots,\xi_n,1,0)\prod_{j=1}^n V(\xi_j)\,\ud \xi_1\cdots
    \ud \xi_n.
  \end{equation}
  Therefore, $v_R(V;1,0)$ depends only on the values $V(\xi)$ for $\xi \in
  [0,1]\times[-R,R]$.}
\end{remark}

\begin{remark}
  {\rm From the chaos expansion \eqref{E:chaos-v}, we deduce that for any
  non-negative bounded potential $V$, we have
  \[
    1\le v_R(V;t,x)\le\cosh\bigl(Rt\sqrt{\|V\|_\infty}\bigr)
  \]
  by~\cite[Lemma~2.1,
  (2.3)]{balan.chen.ea:26:exact*1}; in its proof the mass of
  $f_{n,R}(\cdot,t,x)$ is $t^{2n}/(2n)!$. It follows that there exist some
  constants $c_1,c_2 > 0$ such that for all $t \in [0,1]$ and $x \in \bR$,
  \begin{equation}\label{E:UB-vR}
    v_R(V;t,x)  \le c_1 \exp\big(c_2 R \sqrt{\|V\|_{\infty} t}\big).
  \end{equation}
  }
\end{remark}

We recall Hypothesis~A of~\cite{balan.chen.ea:26:exact*1}.
\begin{hypA}\label{H:hypA}
  $V:[0,1]\times\bR \to [0,\infty)$ is a measurable function which satisfies
  the following conditions:\\
  a) there exists $B > 0$ such that $V(t,x) \le B$ for all
  $(t,x) \in [0,1]\times\bR$;\\
  b) $V(t,\cdot)$ is infinitely differentiable for all
  $t \in [0,1]$;\\
  c) there exists a non-decreasing function $w:[0,1] \to [0,\infty)$, with $w(0)
  = 0$ and $\lim_{h\downarrow0}w(h) = 0$, such that
  \begin{equation}\label{E:time-modulus}
    \|V(t,\cdot)-V(s,\cdot)\|_\infty
    \le w(|t-s|)
    \quad\mbox{for all }t,s\in[0,1].
  \end{equation}
\end{hypA}

The following is Theorem~1.1(i) of~\cite{balan.chen.ea:26:exact*1}.
\begin{theorem}\label{T:log-v-th}
  Suppose that $V$ satisfies Hypothesis A. Then
  \[
    \lim_{R\to\infty}\frac{1}{R}\log v_R(V;1,0)=\Phi(V),
  \]
  where $\Phi(V)$ is given by~\eqref{E:def-Phi-V}.
\end{theorem}

The following is Theorem~1.1(ii) of~\cite{balan.chen.ea:26:exact*1}.
\begin{theorem}\label{T:uniform-upper-bound}
  Fix $B > 0$ and a time modulus $w$ as in Hypothesis~A(c). For every
  $\varepsilon > 0$, there exists $R_0 = R_0(B,w,\varepsilon) > 0$ such that,
  for every potential $V$ satisfying Hypothesis~A with bound $B$ and time
  modulus $w$,
  \[
    \frac{1}{R}\log v_R(V;1,0)\le \Phi(V)+\varepsilon
    \qquad\text{for all } R\ge R_0,
  \]
  where $\Phi(V)$ is given by~\eqref{E:def-Phi-V}.
\end{theorem}

In particular, the limit in Theorem~1.1(i)
of~\cite{balan.chen.ea:26:exact*1}
gives the following lower bound.
\begin{theorem}\label{T:lower-bound}
  Suppose that $V$ satisfies Hypothesis~A. Then
  \[
    \liminf_{R\to\infty}\frac{1}{R}\log v_R(V;1,0)\ge\Phi(V).
  \]
\end{theorem}

\section{Upper bound for the second moment}
\label{S:section-UB}

In this section, we prove the following result.

\begin{theorem}\label{T:UB-th}
  For the constant $D$ given in~\eqref{E:def-D}, we have:
  \[
    \limsup_{R \to \infty} \frac{1}{R}\log \|U_{R}^{Q}(1,0)\|_2 \le D.
  \]
  Consequently, by relation \eqref{E:u-U-sup} with $p = 2$,
  \[
    \limsup_{t \to \infty} t^{-\frac{4-\alpha-\beta}{3-\alpha}} \log
    \|u_{\theta}(t,x)\|_2 \le \theta^{\frac{1}{3-\alpha}}D.
  \]
\end{theorem}

In Section~\ref{SS:section-UB-smooth}, we show that a similar upper bound holds
for the solution $U_R^{Q_N}$ of the equation with smooth noise $\dot{W}^{Q_N}$.
In Section~\ref{SS:section-cutoff} we show how to remove the cut-off, and prove
Theorem~\ref{T:UB-th} by letting $N \to \infty$.

\subsection{Upper bound for equation with smooth noise}\label{SS:section-UB-smooth}

The goal of this section is to prove the following result.

\begin{theorem}\label{T:UB-smooth}
  For any integer $N \geq 2$,
  \[
    \limsup_{R \to \infty} \frac{1}{R}\log \|U_{R}^{Q_N}(1,0)\|_2 \le D_N,
  \]
  where
  \[
    D_N=\sup_{0\le f \in \cK_N}\left\{\Phi(V_{f}^N)
    -\frac{1}{2}\|f\|_{\cK_N}^2 \right\}
  \]
  and
  \begin{equation}\label{E:def-VfN}
    V_{f}^N (\xi)=\int_{E_N}K_{\xi}(e)f(e)\fm(\ud e), \quad \xi \in [0,1]
    \times
    \bR.
  \end{equation}
\end{theorem}

The proof of Theorem~\ref{T:UB-smooth} is based on a sequence of preliminary
results.

\begin{lemma}\label{L:lem-CjkN}
  For any integers $N \geq 2$, $j \geq 0$, $k \geq
  0$ and for any $(t,x) \in [0,1] \times \bR$,
  \begin{equation}\label{E:def-Cjk}
    \|\partial_t^j \partial_x^k K_{t,x}\|_{\cK_N}^2
    =C_{j,k,N}^2<\infty.
  \end{equation}
\end{lemma}

\begin{proof}
  Note that
  \begin{equation}\label{E:deriv-K}
    \partial_t^j \partial_x^k K_{t,x}(a,y,b,z)=\fg_a^{(j)}(t-y)
    \fg_b^{(k)}(x-z),
  \end{equation}
  where $\fg_a^{(j)}(x)\coloneqq \frac{\ud^j}{\ud x^j} \fg_a(x)$. For any $a >
  0$ and $x \in \bR$, $\fg_a(x) = a^{\frac{1}{4}} \fg_1(\sqrt{a} x)$. Hence,
  \begin{equation}\label{E:gaj}
    \fg_a^{(j)}(x)=a^{\frac{1}{4}+\frac{j}{2}} \fg_1^{(j)} (\sqrt{a} x),
  \end{equation}
  and $\int_{\bR}|\fg_a^{(j)}(x)|^2 \ud x = a^j \int_{\bR} |\fg_1^{(j)}(y)|^2\ud
  y =: a^j C_j < \infty$. We have:
  \begin{align*}
    \|\partial_t^j \partial_x^k K_{t,x}\|_{\cK_N}^2& =\int_{E_N}(
    \partial_t^j
    \partial_x^k K_{t,x}(a,y,b,z) )^2 w_{\beta}(a) w_{\alpha}(b) \ud a
    \ud y\ud b\ud z\\
    &=\int_{\frac{1}{N}}^{N} \int_{\frac{1}{N}}^{N} \left(\int_{\bR}
    |\fg_a^{(j)}(t-y)|^2\ud y\right)
    \left(\int_{\bR} |\fg_b^{(k)}(x-z)|^2\ud z\right) w_{\beta}(a)
    w_{\alpha}(b) \ud a\ud b \\
    &=C_j C_k \left(\int_{\frac{1}{N}}^{N} a^j w_{\beta}(a) \ud a \right)
    \left(\int_{\frac{1}{N}}^{N} b^k w_{\alpha}(b)
    \ud b\right)=:C_{j,k,N}^2<\infty.
  \end{align*}
  This proves Lemma~\ref{L:lem-CjkN}.
\end{proof}

Recall the definition \eqref{E:def-WQN} of $\dot{W}^{Q_N}$.

\begin{lemma}\label{L:deriv-XN}
  For any integer $N \geq 2$, the map $(t,x) \mapsto \dot{W}^{Q_N}(t,x) =
  B_N(K_{t,x}) \in L^2(\Omega)$ is smooth. Moreover, for any integers $j \geq
  0$, $k \geq 0$ and for any $(t,x) \in [0,1] \times \bR$,
  \[
    \partial_t^j \partial_x^k \dot{W}^{Q_N}(t,x)=B_N(\partial_t^j
    \partial_x^k
    K_{t,x}) \quad \mbox{in} \quad L^2(\Omega).
  \]

\end{lemma}

\begin{proof}
  Fix $N\geq2$. The formula for $K_{t,x}$ extends to $(t,x) \in \bR^2$, and
  Lemma~\ref{L:lem-CjkN} remains valid there. In Lemma~\ref{L:lem-CjkN}, the
  derivatives of $K_{t,x}(e)$ were taken for each fixed $e$. Here we first
  justify differentiation in the $\cK_N$ norm.

  For any integers $j,k\geq0$, the fundamental theorem of calculus at each
  fixed $e$ gives
  \[
    \partial_t^j\partial_x^kK_{t+h,x}
    -\partial_t^j\partial_x^kK_{t,x}
    =h\int_0^1\partial_t^{j+1}\partial_x^kK_{t+sh,x}\,\ud s.
  \]
  By Lemma~\ref{L:lem-CjkN}, we have
  \begin{equation*}
    \begin{aligned}
      \|\partial_t^j\partial_x^kK_{t+h,x}
      -\partial_t^j\partial_x^kK_{t,x}\|_{\cK_N}^2&=h^2\int_{E_N}
      \left|\int_0^1
      \partial_t^{j+1}\partial_x^kK_{t+sh,x}(e)\,\ud s
      \right|^2\fm(\ud e)\\
      &\le h^2\int_{E_N}\int_0^1
      |\partial_t^{j+1}\partial_x^kK_{t+sh,x}(e)|^2
      \,\ud s\,\fm(\ud e)\\
      &=h^2\int_0^1
      \|\partial_t^{j+1}\partial_x^kK_{t+sh,x}\|_{\cK_N}^2\,
      \ud s=h^2C_{j+1,k,N}^2.
    \end{aligned}
  \end{equation*}
  Similarly,
  \[
    \|\partial_t^j\partial_x^kK_{t,x+h}
    -\partial_t^j\partial_x^kK_{t,x}\|_{\cK_N}^2
    \le h^2C_{j,k+1,N}^2.
  \]
  Consequently,
  \[
    \|\partial_t^j\partial_x^kK_{t+h,x+\ell}
    -\partial_t^j\partial_x^kK_{t,x}\|_{\cK_N}
    \le |h|C_{j+1,k,N}
    +|\ell|C_{j,k+1,N}
    \longrightarrow0
  \]
  as $(h,\ell) \to (0,0)$. Thus, $(t,x)\mapsto\partial_t^j\partial_x^kK_{t,x}$
  is continuous in $\cK_N$.

  We now verify that these pointwise derivatives are also derivatives in the
  $\cK_N$ norm. Using the same integral identity and applying the first estimate
  with $j+1$ in place of $j$, we get
  \[
    \left\|
    \frac{\partial_t^j\partial_x^kK_{t+h,x}
    -\partial_t^j\partial_x^kK_{t,x}}{h}
    -\partial_t^{j+1}\partial_x^kK_{t,x}
    \right\|_{\cK_N}^2
    \le \frac{h^2}{3}C_{j+2,k,N}^2
    \longrightarrow0.
  \]
  The same argument applies to differentiation with respect to $x$. Hence
  $(t,x)\mapsto K_{t,x} \in \cK_N$ is smooth.

  Finally, since $B_N$ is a linear isometry, the formula in the statement
  follows, and these derivatives are continuous. This completes the proof of
  Lemma~\ref{L:deriv-XN}.
\end{proof}

We denote $C_R = [0,1] \times [-R,R]$. For any function $f:S \to \bR$, we
denote
\[
  \|f\|_{\infty, S}=\sup_{x \in S} |f(x)|.
\]

\begin{lemma}\label{L:4-mom-Z}
  Recall definition~\eqref{E:def-ZN} of $Z_N$. For any integer $N \geq 2$ and
  for any $R > 0$,
  \[
    \bE\|Z_N\|_{\infty,C_R}^4 \le L_N(R+1),
  \]
  where $L_N = 256 \sum_{j,k \in \{0,1\}} C_{j,k,N}^4$ and $C_{j,k,N}$ is given
  by~\eqref{E:def-Cjk}.
\end{lemma}

\begin{proof}
  Using~\eqref{E:def-WQN}, we can assume that
  \[
    Z_N(\xi)=\frac{B_{N,1}(K_{\xi})+ iB_{N,2}(K_{\xi}) }{\sqrt{2}}, \quad
    \mbox{for all} \ \xi=(t,x) \in [0,\infty) \times \bR,
  \]
  where $B_{N,1}$ and $B_{N,2}$ are independent copies of $B_N$. By
  Lemma~\ref{L:deriv-XN},
  \[
    \partial_t^j \partial_x^k Z_N(\xi)=\frac{B_{N,1}(\partial_t^j
    \partial_x^k
    K_{\xi}) + i B_{N,2}(\partial_t^j \partial_x^k K_{\xi}) }{\sqrt{2}}.
  \]
  Hence,
  \begin{equation}
    \bE|\partial_t^j \partial_x^k Z_N(\xi)|^4
    = \frac12 \bE|B_{N,1}(\partial_t^j
    \partial_x^k K_\xi)|^4
    \label{E:partial-Z}
    + \frac12
    \left(\bE|B_{N,1}(\partial_t^j
    \partial_x^k K_\xi)|^2\right)^2
    = 2C_{j,k,N}^4,
  \end{equation}
  using Lemma~\ref{L:lem-CjkN}, and the fact that $B_N(\partial_t^j
  \partial_x^k
  K_{\xi})$ has $N(0,C_{j,k,N}^2)$ distribution.

  Apply the same Gaussian increment argument to each derivative in the rectangle
  identity, using the next derivatives in Lemma~\ref{L:deriv-XN}. The integrated
  identity holds in $L^2(\Omega)$ at each point by the fundamental theorem of
  calculus, hence almost surely on a countable dense set and then everywhere for
  the continuous versions. Let $D_j = [0,1] \times [j,j+1]$. By
  Lemma~\ref{L:deriv-XN}, for any integer $j \geq 0$, and for $(t,x) \in D_j$,
  \[
    Z_N(t,x)=Z_N(0,j)+\int_j^x \partial_x Z_N(0,v)\ud v+\int_0^t \partial_t
    Z_N(u,j)\ud u+\int_0^t \int_j^x \partial_t \partial_x Z_N(u,v)\ud u\ud v.
  \]
  Using H\"older's inequality, followed by~\eqref{E:partial-Z}, we get:
  \begin{align}
    \nonumber
    \bE\|Z_N\|_{\infty,D_j}^{4} & \le 4^3 \left\{
    \bE|Z_N(0,j)|^4+\int_j^{j+1}
    \bE|\partial_x Z_N(0,v)|^4 \ud v+\int_0^1 \bE|\partial_t Z_N(u,j)|^4
    \ud u
    \right. \\
    \nonumber
    & \quad \quad \left. + \int_0^1 \int_{j}^{j+1} \bE|\partial_t \partial_x
    Z_N(u,v)|^4\ud u\ud v \right\}\\
    \label{E:bound-EZ}
    & = 128\big(C_{0,0,N}^4+ C_{0,1,N}^4+C_{1,0,N}^4+C_{1,1,N}^4\big).
  \end{align}
  A similar relation holds for any integer $j < 0$. Note that
  \[
    \|Z_N\|_{\infty,C_R}^4 \le \max_{-\lceil R\rceil\le j<\lceil
    R\rceil} \|Z_N\|_{\infty,D_j}^4 \le \sum_{j=-\lceil
    R\rceil}^{\lceil R\rceil-1} \|Z_N\|_{\infty,D_j}^4.
  \]
  There are $2\lceil R\rceil\le2(R+1)$ cells. The conclusion follows after
  taking expectation, and using~\eqref{E:bound-EZ}. This completes the proof of
  Lemma~\ref{L:4-mom-Z}.
\end{proof}

\begin{lemma}\label{L:lem-12}
  For any integer $N \geq 2$ and for any $R > 0$,
  \[
    \bP \left(\left\| \frac{Z_N}{\sqrt{R}} \right\|_{\infty,C_R} \le
    \varepsilon_R \right)\geq \frac{1}{2} \quad
    {\rm where} \quad \varepsilon_R=\frac{[2L_N(R+1)]^{1/4}}{\sqrt{R}},
  \]
  and $L_N$ is given in Lemma~\ref{L:4-mom-Z}.
\end{lemma}

\begin{proof}
  By Markov's inequality and Lemma~\ref{L:4-mom-Z},
  \[
    \bP \left(\left\| \frac{Z_N}{\sqrt{R}} \right\|_{\infty,C_R} \geq
    \varepsilon_R \right) \le \frac{1}{\varepsilon_R^4} \bE\left\|
    \frac{Z_N}{\sqrt{R}}\right\|_{\infty,C_R}^4 \le
    \frac{1}{\varepsilon_R^4}
    \frac{L_N(R+1)}{R^2}=\frac{1}{2}.
  \]
  This proves Lemma~\ref{L:lem-12}.
\end{proof}

We recall from Appendix~\ref{S:appB} that $\cH_R$ is the Cameron-Martin space of
the process
\[
  Z_R^*(\xi)=\frac{Z_N(\xi)}{\sqrt{R}}, \ \xi\in C_R,
\]
when it is viewed as a Gaussian random element in the Banach space
$\mathbf{B} = C(C_R;\bC)$. Let
\begin{equation}\label{E:def-rho}
  \rho_R(z)=\inf\{\|f\|_{\cK_N^{\bC}}; \|z-T_R f\|_{\infty,C_R} \le
  \varepsilon_{R} \} \quad \mbox{for} \quad z \in \mathbf{B},
\end{equation}
with the convention $\inf \emptyset = \infty$, where $\cK_N^{\bC}$ and $T_R$ are
defined in Appendix~\ref{S:appB}.

\medskip

\emph{Notation.} Since the cutoff $N$
is fixed, we omit it from the notation for $\varepsilon_R$, $Z_R^*$, $\cH_R$,
$T_R$, and $\rho_R$.

\begin{lemma}\label{L:inf-attained}
  If $z \in \mathbf{B}$ is such that $\rho_R(z) < \infty$, then the
  infimum in \eqref{E:def-rho} is attained, i.e. there exists $f
  \in \cK_{N}^{\bC}$ such that
  $\|z-T_R f\|_{\infty,C_R} \le \varepsilon_{R}$ and $\rho_R(z) =
  \|f\|_{\cK_N^{\bC}}$.
\end{lemma}

\begin{proof}
  There exists a sequence $(f_n)_{n \geq 1}$ in $\cK_N^{\bC}$ such that $\|z-T_R
  f_n\|_{\infty,C_R} \le \varepsilon_{R}$ and
  \[
    \rho_R(z) \le \|f_n\|_{\cK_N^{\bC}} \le \rho_R(z)+\frac{1}{n}.
  \]
  Since $(f_n)_{n \geq 1}$ is bounded in $\cK_N^{\bC}$, there exists a
  subsequence $N'\subset \bN$ and $f \in \cK_{N}^{\bC}$ such that
  $f_n \to f$ weakly in $\cK_{N}^{\bC}$, as $n \to \infty$, $n \in N'$, i.e.
  \[
    \lim_{n \to \infty,n \in N'}\langle f_n,g\rangle_{\cK_N^{\bC}} = \langle
    f,g\rangle_{\cK_N^{\bC}} \quad \mbox{for all} \ g \in \cK_{N}^{\bC}.
  \]
  In the real Hilbert space $\cK_N^\bC$, pairing $f_n$ with
  $K_\xi$ and $iK_\xi$ gives the real and imaginary parts of
  $V_{f_n}^N(\xi)$, respectively. Weak convergence against both vectors gives
  \[
    \lim_{n \to \infty,n \in N'}V_{f_n}^N(\xi)=V_f^N(\xi) \quad \mbox{for
    any}
    \ \xi \in [0,1] \times \bR.
  \]

  From Appendix~\ref{S:appB}, we recall that $T_R f$ is the restriction of
  $V_f^N$ to $C_R$. Hence,
  \[
    |z(\xi)-V_f^N(\xi)| \le |z(\xi)-(T_R
    f_n)(\xi)|+|V_{f_n}^N(\xi)-V_f^N(\xi)|\le
    \varepsilon_R+|V_{f_n}^N(\xi)-V_f^N(\xi)|,
  \]
  for any $\xi \in C_R$. Letting $n \to \infty$ (for $n \in N'$), we get
  $\|z-V_f^N\|_{\infty,C_R} \le \varepsilon_R$. By the definition of
  $\rho_R(z)$, it follows that $\rho_R(z) \le \|f\|_{\cK_N^{\bC}}$.

  On the other hand, since $\|\cdot\|_{\cK_N^{\bC}}$ is lower semi-continuous
  for the weak topology,
  \[
    \rho_R(z)=\liminf_{n \to \infty,n \in N'}\|f_n\|_{\cK_N^{\bC}} \geq
    \|f\|_{\cK_N^{\bC}}.
  \]
  Hence, $\rho_R(z) = \|f\|_{\cK_N^{\bC}}$.
  This proves Lemma~\ref{L:inf-attained}.
\end{proof}

The same argument shows that $\rho_R$ is lower semicontinuous. Indeed, if
$z_n\to z$ uniformly and $\rho_R(z_n)\le a<\infty$, choose minimizing
controls of norm at most $a$. Along a weakly convergent subsequence their
images converge pointwise, so the limit control has norm at most $a$ and
approximates $z$ within $\varepsilon_R$ at every point. Thus every sublevel
set is closed and $\rho_R$ is Borel measurable.

We will use the following result (Borell-Sudakov-Cirel'son Theorem).

\begin{theorem}[Theorem 4.3 of \cite{ledoux:96:isoperimetry}]
  \label{T:BSC-th}
  Let $Z$ be a zero-mean Gaussian random element in a real separable Banach
  space $B$. Let $B_{\varepsilon}$ be the closed ball of radius $\varepsilon >
  0$ centered at origin in the Cameron-Martin space of $Z$, and $A \subseteq B$
  be a closed set. Then
  \[
    \bP(Z \in A+B_{\varepsilon}) \geq \Phi\big(\Phi^{-1} (\bP(Z \in
    A))+\varepsilon\big),
  \]
  where $\Phi$ is the standard normal cumulative distribution function, with
  $\Phi^{-1}(0)=-\infty$ and $\Phi^{-1}(1)=+\infty$, interpreting
  $\Phi(\pm\infty)$ by continuity.
\end{theorem}

\begin{lemma}\label{L:rho-tail}
  For any $a > 0$,
  \[
    \bP\big(\rho_R(Z_R^*) \geq a\big) \le e^{-R a^2}.
  \]
  Consequently, $\bP\big(\rho_R(Z_R^*) < \infty\big) = 1$.
\end{lemma}

\begin{proof}
  Let $z \in \mathbf{B}$ be arbitrary. Note that $\rho_R(z) < a$ if
  and only if there exists $f \in \cK_N^{\bC}$ such that
  $\|z-T_Rf\|_{\infty,C_R} \le \varepsilon_R$ and $\|f\|_{\cK_N^{\bC}} < a$,
  i.e. $z = b + h$, where
  \[
    b=z-T_Rf\in A_R \quad \mbox{with} \quad A_R=\{z \in
    \mathbf{B};
    \|z\|_{\infty,C_R} \le \varepsilon_R\},
  \]
  while $h = T_R f$ belongs to $\sqrt{2R}a B_{\cH_R}$, the open ball of radius
  $\sqrt{2R}a$ in $\cH_R$, according to the representation of this ball given by
  Corollary~\ref{C:CM-ball}. (We recall from Appendix~\ref{S:appB} that
  $B_{\cH_R}$ is the open unit ball in $\cH_R$.) In other words,
  \[
    \rho_R(z)<a \quad \mbox{if and only if} \quad z \in A_R+\sqrt{2R}a
    B_{\cH_R}.
  \]

  The set $A_R$ is closed in $\mathbf{B}$, and for any $r \in (0,\sqrt{2R}a)$,
  the closed ball $B_r$ of radius $r$ in $\cH_R$ is contained in $\sqrt{2R} a
  B_{\cH_R}$. Therefore, by applying Theorem~\ref{T:BSC-th}, we obtain that:
  \begin{align*}
    \bP\big(\rho_R(Z_R^*)<a\big)
    & = \bP\big( Z_R^* \in A_R +\sqrt{2R} a B_{\cH_R}\big)       \\
    & \geq \bP\big( Z_R^* \in A_R + B_r\big)
    \geq \Phi \Big(\Phi^{-1} \big( \bP(Z_R^* \in A_R) \big) + r\Big).
  \end{align*}
  Letting $r \uparrow \sqrt{2R}a$ and using the continuity of $\Phi$, we get
  \[
    \bP\big(\rho_R(Z_R^*)<a\big) \geq \Phi \Big(\Phi^{-1} \big( \bP(Z_R^* \in
    A_R) \big) +\sqrt{2R} a\Big).
  \]

  By Lemma~\ref{L:lem-12}, $\Phi^{-1} \big( \bP(Z_R^* \in A_R)
  \big) \geq \Phi^{-1}(\frac{1}{2}) = 0$, since $\Phi^{-1}$ is
  non-decreasing. Since $\Phi$ is also non-decreasing, we obtain that
  $\bP\big(\rho_R(Z_R^*) < a\big) \geq \Phi(\sqrt{2R}a)$, and hence
  \[
    \bP\big(\rho_R(Z_R^*)\geq a\big) \le 1-\Phi(\sqrt{2R}a)\le
    e^{-R
    a^2}.
  \]
  This proves Lemma~\ref{L:rho-tail}.
\end{proof}

\begin{lemma}\label{L:Z-V-lem}
  Let $\Omega_{R} = \{\rho_R(Z_R^*) < \infty\}$. For any $\omega \in
  \Omega_R$,
  there exists a function $f_R(\omega) \in \cK_N^{\bC}$ such that on the event
  $\Omega_R$, $\rho_R(Z_R^*) = \|f_R\|_{\cK_N^\bC}$ and
  \[
    |Z_R^*(\xi)| \le V_{|f_R|}^N(\xi)+\varepsilon_R \quad \mbox{for all}
    \ \xi
    \in C_R.
  \]
  ($\Omega_R$ and $f_R$ depend also on $N$.)
\end{lemma}

\begin{proof}
  By Lemma~\ref{L:inf-attained}, for any $\omega \in \Omega_R$, there exists a
  function $f_R(\omega) \in \cK_{N}^{\bC}$ such that
  \[
    \|Z_R^*(\omega)-T_R f_R(\omega)\|_{\infty,C_R} \le \varepsilon_{R}
    \quad
    \mbox{and} \quad \rho_R(Z_R^*(\omega))=\|f_R(\omega)\|_{\cK_N^{\bC}}.
  \]
  Recall that $T_R f_R(\xi) = V_{f_R}^N (\xi)$ for any $\xi \in C_R$. Hence, on
  the event $\Omega_R$, for any $\xi \in C_R$,
  \[
    |Z_R^*(\xi)| \le |Z_R^*(\xi)-V_{f_R}^N(\xi)|+|V_{f_R}^N(\xi)| \le
    \varepsilon_R+ V_{|f_R|}^N(\xi).
  \]
  This proves Lemma~\ref{L:Z-V-lem}.
\end{proof}

\begin{corollary}\label{C:vv-cor}
  On the event $\Omega_R$,
  \[
    \left|v_R\left(\frac{Z_N}{\sqrt{R}};1,0\right)\right| \le v_R
    (V_{|f_R|}^N+\varepsilon_R;1,0).
  \]
\end{corollary}

\begin{proof}
  This follows by relation \eqref{E:int-CR} and
  Lemma~\ref{L:Z-V-lem}.
  In~\eqref{E:int-CR} all factors $f_{n,R}$ are nonnegative and
  all points lie in $C_R$; taking absolute values term by term permits the bound
  of Lemma~\ref{L:Z-V-lem}.
  This proves Corollary~\ref{C:vv-cor}.
\end{proof}

\begin{lemma}\label{L:VN-HypA}
  For any non-negative function $f \in \cK_N$,
  $V_f^N$ satisfies Hypothesis~\ref{H:hypA}.
\end{lemma}

\begin{proof}
  Clearly, $V_f^N$ is non-negative. By the Cauchy--Schwarz
  inequality and Lemma~\ref{L:lem-CjkN},
  \begin{equation}\label{E:Vf-bded}
    V_f^N(t,x) \le C_{0,0,N}\|f\|_{\cK_N}.
  \end{equation}
  This proves part a) of Hypothesis~\ref{H:hypA}.
  By the proof of Lemma~\ref{L:deriv-XN}, the map
  $(t,x)\mapsto K_{t,x}\in\cK_N$ is smooth. Since
  $V_f^N(t,x)=\langle K_{t,x},f\rangle_{\cK_N}$,
  it follows that $V_f^N$ is smooth and, by the
  Cauchy--Schwarz inequality,
  \[
    |\partial_t^j\partial_x^k V_f^N(t,x)|
    \le C_{j,k,N}\|f\|_{\cK_N}.
  \]
  This proves part b). Finally, by the mean value theorem,
  \[
    \|V_f^N(t,\cdot)-V_f^N(s,\cdot)\|_\infty
    \le C_{1,0,N}\|f\|_{\cK_N}|t-s|,
  \]
  which proves part c).
  This proves Lemma~\ref{L:VN-HypA}.
\end{proof}

\begin{lemma}\label{L:delta-lem}
  For any $\delta > 0$, there exist $\eta \in (0,1)$ and $R_0 > 0$ such that
  \[
    \frac{1}{R} \log v_R (V_{f}^N+\varepsilon_R;1,0) \le D_N+
    \frac{1-\eta}{2}\|f\|_{\cK_N}^2+\delta,
  \]
  for any  $R > R_0$ and for any non-negative function $f \in \cK_N$.
\end{lemma}

\begin{proof}
  Let $K = 4(C_{0,0,N}+1)$, where $C_{0,0,N}$ is given by~\eqref{E:def-Cjk}.
  Note that $\varepsilon_R \le 1$ for $R$ large enough.
  Let $f \in \cK_N$ be a non-negative function.
  By Lemma~\ref{L:VN-HypA}, $V_f^N$ satisfies
  Hypothesis~\ref{H:hypA}.

  For $\|f\|_{\cK_N} \le K$ and $\varepsilon_R\le1$, the preceding estimates give bound
  $C_{0,0,N}K+1$ and time modulus $h\mapsto C_{1,0,N}Kh$ for
  $V_f^N+\varepsilon_R$.
  By Theorem~\ref{T:uniform-upper-bound}, one threshold $R_0$ works
  for the whole ball. Therefore, for any $\delta > 0$ there exists $R_1 > 0$
  (depending on $\delta$), such that for all $\|f\|_{\cK_N}\le K$,
  \[
    \frac{1}{R}\log v_R\big(V_f^N+\varepsilon_R;1,0\big) \le
    \Phi(V_f^N+\varepsilon_R)+\frac{\delta}{4} \quad \mbox{for all} \ R>R_1.
  \]
  Note that $\Phi(V+\varepsilon) \le \Phi(V)+\sqrt{\varepsilon}$
  by
  shifting each Rayleigh supremum,
  for any non-negative function $V$ and for any $\varepsilon > 0$.
  We choose $R_2 > R_1$ (depending on $\delta$) such that
  $\varepsilon_R\le1$ and $\sqrt{\varepsilon_R} \le
  \frac{\delta}{4}$ for any $R > R_2$. Then
  \[
    \frac{1}{R}\log v_R\big(V_f^N+\varepsilon_R;1,0\big) \le
    \Phi(V_f^N)+\frac{\delta}{2}\le
    D_N+\frac{1}{2}\|f\|_{\cK_N}^2+\frac{\delta}{2} \quad \mbox{for all}
    \ R>R_2.
  \]

  We consider two cases.

  \medskip

  {\bf Case 1.} $\|f\|_{\cK_N} \le K$. In this case, for any $\eta <
  \frac{\delta}{K^2}$,
  \[
    \frac{1}{2}\|f\|_{\cK_N}^2+\frac{\delta}{2} \le
    \frac{1-\eta}{2}\|f\|_{\cK_N}^2 +\frac{\eta K^2}{2}+\frac{\delta}{2}\le
    \frac{1-\eta}{2}\|f\|_{\cK_N}^2+\delta.
  \]

  {\bf Case 2.} $\|f\|_{\cK_N} > K$. In this case, by~\eqref{E:Vf-bded},
  \[
    V_f^N(\xi)+\varepsilon_R \le C_{0,0,N}\|f\|_{\cK_N}+1 \le
    \frac{\|f\|_{\cK_N}^2}{4}+1 \le \frac{\|f\|_{\cK_N}^2}{2} \le
    \frac{\|f\|_{\cK_N}^4}{16}.
  \]

  From~\eqref{E:UB-vR}, there exist some constants $c_1,c_2 >
  0$
  ($c_1 = c_2 = 1$ is possible since $\cosh(Rt\sqrt M) \le
  e^{R\sqrt{Mt}}$ for $M\ge0$, $0 \le t\le1$) such that
  \[
    v_R\big(V_f^N+\varepsilon_R;1,0\big) \le c_1
    \exp(c_2 R \sqrt{\|V_f^N\|_{\infty}+\varepsilon_R}).
  \]
  Then for any $R > 0$ with $\varepsilon_R\le1$,
  \[
    \frac{1}{R} \log v_R\big(V_f^N+\varepsilon_R;1,0\big) \le
    \frac{\log c_1}{R} + c_2 \frac{\|f\|_{\cK_N}^2}{4}.
  \]

  Choose $0 < \eta < \min(1/2,\delta/K^2)$. With
  $c_1 = c_2 = 1$, the last bound is at most
  $(1-\eta)\|f\|_{\cK_N}^2/2$. The zero control gives $D_N\ge0$,
  which proves the assertion also in Case~2. Take $R_0 = R_2$. This proves
  Lemma~\ref{L:delta-lem}.
\end{proof}

Now we are ready to prove Theorem~\ref{T:UB-smooth}.

\begin{proof}[Proof of Theorem~\ref{T:UB-smooth}]
  Let $\delta > 0$ be arbitrary, and $\eta$ and $R_0$ be given by
  Lemma~\ref{L:delta-lem}. By Corollary~\ref{C:vv-cor}, Lemma~\ref{L:delta-lem}
  and Lemma~\ref{L:Z-V-lem}, for any $R > R_0$, on the event $\Omega_R$,
  \begin{align*}
    \left|v_R\left( \frac{Z_N}{\sqrt{R}};1,0\right)\right|
    \le v_R (V_{|f_R|}^N+\varepsilon_R;1,0)
    \le \exp \left\{ R\Big(D_N+\delta+\frac{1-\eta}{2} \|f_R\|_{\cK_N^\bC}^2 \Big) \right\} \\
    = \exp \left\{ R\left[D_N+\delta+\frac{1-\eta}{2}
    \rho_R^2\Big(\frac{Z_N}{\sqrt{R}}\Big) \right] \right\}.
  \end{align*}
  Combining this estimate with Proposition~\ref{P:uQ-vR}, we infer that for any
  $R > R_0$,
  \[
    \bE|U_R^{Q_N}(1,0)|^2=\bE\Big|v_R\left(
    \frac{Z_N}{\sqrt{R}};1,0\right)\Big
    |^2 \le e^{2R (D_N+\delta)}\bE\big[e^{(1-\eta)R
    \rho_R^2\big(\frac{Z_N}{\sqrt{R}}\big)}\big].
  \]
  To compute the last expectation, we use the fact that:
  \[
    \bE[e^{cX^2}] =1+\int_0^{\infty} 2ca e^{ca^2} P(X \geq a)\ud a,
  \]
  for any non-negative random variable $X$ and for any $c > 0$. Using
  Lemma~\ref{L:rho-tail}, we get:
  \begin{align*}
    \bE\big[e^{(1-\eta)R \rho_R^2\big(\frac{Z_N}{\sqrt{R}}\big)}\big] \le 1
    +2(1-\eta) R \int_0^{\infty} a e^{-\eta R a^2}\ud a=\frac{1}{\eta}.
  \end{align*}
  It follows that for any $R > R_0$, $\|U_R^{Q_N}(1,0)\|_2 \le e^{R(D_N+\delta)}
  \eta^{-1/2}$, and hence,
  \[
    \frac{1}{R}\log \|U_R^{Q_N}(1,0)\|_2 \le D_N+\delta-\frac{1}{2R} \log
    \eta.
  \]
  We take the limsup as $R \to \infty$, then we let $\delta \to 0$. This proves
  Theorem~\ref{T:UB-smooth}.
\end{proof}

\subsection{Removing the cut-off}\label{SS:section-cutoff}

In this section, we give the proof of Theorem~\ref{T:UB-th}.

\medskip

We begin by presenting some properties of the second-moment functional
\[
  M_R(Q)=\bE|U_R^{Q}(1,0)|^2
\]
associated with the solution $U_R^Q$ of equation~\eqref{E:HAM-R}, which will be
used for passing from $Q_N$ to $Q$ in the proof of Theorem~\ref{T:UB-th}.

\medskip

For a nonnegative, nonnegative definite covariance $\rho$, let
$C_{n,R}=n!\widetilde f_{n,R}(\cdot,1,0)$ and define
\[
  M_R(\rho)=1+\sum_{n\ge1}\frac{R^{3n}}{n!}
  \int C_{n,R}(\pmb{\xi})C_{n,R}(\pmb{\eta})
  \prod_{j=1}^n\rho(\xi_j-\eta_j)\,\ud\pmb{\xi}\,\ud\pmb{\eta}.
\]
Here, as stated at the end of Section~\ref{SS:motivation}, the integral is over
$[0,\infty)\times\bR$ in each of the variables $\xi_j$ and $\eta_j$. This
nonnegative series may be infinite. When the square-integrable solution exists,
the chaos isometry gives $M_R(\rho)=\bE|U_R^\rho(1,0)|^2$. The series is
monotone in $\rho$.

The following result shows that $M_{R}(\cdot)$ is log-sub-additive.

\begin{lemma}\label{L:log-subad}
  For any non-negative and non-negative definite functions $\rho_1$ and
  $\rho_2$,
  \[
    M_{R}(\rho_1+\rho_2)\le M_R(\rho_1)M_{R}(\rho_2).
  \]
\end{lemma}

\begin{proof}
  Let $C_{n,R} = n! \widetilde{f}_{n,R}(\cdot,1,0)$. It can be proved that
  \begin{equation}\label{E:Cnm}
    C_{n+m,R}(\pmb{z},\pmb{w})\le
    C_{n,R}(\pmb{z})C_{m,R}(\pmb{w}),
  \end{equation}
  for any $\pmb{z} = (z_1,\ldots,z_n) \in ([0,\infty) \times \bR)^n$ and
  $\pmb{w} = (w_1,\ldots,w_m) \in ([0,\infty) \times \bR)^m$. In chronological
  order the kernel is $(2R)^{-\ell}$ for a chain of $\ell$ points with all cone
  constraints satisfied, including the edge to $(1,0)$, and is zero otherwise.
  Only this order contributes. On deleting points, the triangle inequality
  preserves each resulting edge constraint. Apply this once to retain $\pmb z$
  and once to retain $\pmb w$. The weights then multiply to $(2R)^{-(n+m)}$, as
  required. A tie gives zero on the left by strict time ordering; $C_{0,R} =
  1$.

  Note that
  \begin{align*}
    M_R(\rho_1+\rho_2) = 1+\sum_{n\geq 1}\frac{R^{3n}}{n!} \int C_{n,R}(\pmb{\xi}) C_{n,R}(\pmb{\eta})
    \prod_{j=1}^{n}
    \big(\rho_1(\xi_j-\eta_j)+ \rho_2(\xi_j-\eta_j)\big) \ud \pmb{\xi}
    \ud \pmb{\eta}.
  \end{align*}
  To evaluate the product, we use the formula:
  \[
    \prod_{j=1}^n (a_j+b_j)=\sum_{I \subseteq \{1,\ldots,n\}} \prod_{j\in
    I}a_j
    \prod_{j\in I^c}b_j.
  \]

  For the term corresponding to the set $I = \{i_1,\ldots,i_k\}$ with $I^c =
  \{j_1,\ldots,j_{n-k}\}$, we use the change of variables: $\xi_{i_{p}} =
  z_{p}$, $\eta_{i_{p}} = z_{p}'$ for $p = 1,\ldots,k$, and $\xi_{j_q} = w_q$,
  $\eta_{j_q} = w_q'$ for $q = 1,\ldots,n-k$. We denote $\pmb{z} =
  (z_1,\ldots,z_k)$, $\pmb{w} = (w_1,\ldots,w_{n-k})$, and similarly $\pmb{z}'$,
  $\pmb{w}'$. This change of variables shows that all the $\binom{n}{k}$
  integrals are the same. Hence,
  \begin{align*}
    M_R(\rho_1+\rho_2) &
    =\sum_{n\geq 0}\frac{R^{3n}}{n!} \sum_{k=0}^{n}\binom{n}{k}
    \int C_{n,R}(\pmb{z},\pmb{w})
    C_{n,R}(\pmb{z}',\pmb{w}') \\
    & \quad \quad \quad \quad \times \prod_{p=1}^{k}
    \rho_1(z_{p}-z'_{p})\prod_{q=1}^{n-k}\rho_2(w_q-w_q')
    \ud \pmb{z}
    \ud \pmb{w} \ud \pmb{z}'\ud \pmb{w}' \\
    & =\sum_{k\geq 0} \sum_{m\geq 0} R^{3(k+m)} \frac{1}{k!
    m!}\int C_{k+m,R}(\pmb{z},\pmb{w})
    C_{k+m,R}(\pmb{z}',\pmb{w}') \\
    & \quad \quad \quad \quad \times \prod_{p=1}^{k}
    \rho_1(z_{p}-z'_{p})\prod_{q=1}^{m}\rho_2(w_{q}-w_{q}')
    \ud \pmb{z} \ud \pmb{w} \ud \pmb{z}'\ud \pmb{w}'.
  \end{align*}
  Using inequality \eqref{E:Cnm}, we have $C_{k+m,R}(\pmb{z},\pmb{w}) \le
  C_{k,R}(\pmb{z})C_{m,R}(\pmb{w})$. Hence,
  \begin{align*}
    M_R(\rho_1+\rho_2) & \le \left(\sum_{k\geq 0} R^{3k}
    \frac{1}{k!}\int
    C_{k,R}(\pmb{z})C_{k,R}(\pmb{z}')
    \prod_{p=1}^{k}
    \rho_1(z_{p}-z'_{p}) \ud \pmb{z} \ud \pmb{z}' \right) \\
    & \quad \times \left(\sum_{m\geq 0} R^{3m} \frac{1}{m!}\int C_{m,R}(\pmb{w})C_{m,R}(\pmb{w}')
    \prod_{q=1}^{m}
    \rho_2(w_{q}-w'_{q}) \ud \pmb{w} \ud \pmb{w}' \right)
    \\ &= M_R(\rho_1)M_R(\rho_2).
  \end{align*}
  This proves Lemma~\ref{L:log-subad}.
\end{proof}

Note that, by~\cite[(70)]{balan.conus:16:intermittency}, for any $a > 0$, there
exist some constants $c_{a,1} > 0$ and $c_{a,2} > 0$ such that
\begin{equation}\label{E:series-x}
  \sum_{n\geq 0}\frac{x^n}{\Gamma(an+1)} \le c_{a,1} \exp(c_{a,2} x^{1/a})
  \quad \mbox{for all} \ x>0.
\end{equation}

The next result shows that for our covariance $Q$ given by~\eqref{E:def-Q},
$M_R(Q)$ has exponential growth in $R$.

\begin{lemma}\label{L:MR-expo}
  Let $Q_{a,b}(t,x) = |t|^{-b}|x|^{-a}$ for some $a,b \in [0,1)$. For any $c >
  0$,
  \begin{equation}\label{E:bound-MQ}
    M_R(cQ_{a,b}) \le \sum_{n\geq 0}\frac{\big(cA_{a,b}R^{3-a}\big)^n
    }{\Gamma((3-a)n+1)} \le c_{3-a,1} \exp\left\{ c_{3-a,2}
    \big(cA_{a,b}\big)^{\frac{1}{3-a}} R\right\},
  \end{equation}
  where $A_{a,b} > 0$ depends on $(a,b)$, and $c_{3-a,1},c_{3-a,2} > 0$ depend
  on $a$.
\end{lemma}

\begin{proof}
  We treat only the case $b \in (0,1)$, $b = 0$ being similar. Let $\cH_{a,b} =
  \cH^{Q_{a,b}}$. Then
  \[
    M_R(cQ_{a,b})=1+\sum_{n\geq 1}R^{3n} c^n n!
    \|\widetilde{f}_{n,R}(\cdot,1,0)\|_{\cH_{a,b}^{\otimes n}}^2.
  \]

  Let $\cH_a$ be the completion of $C_c^{\infty}(\bR)$ with respect to the inner
  product:
  \[
    \langle \varphi,\psi \rangle_{\cH_a}=\int_{\bR^2}\varphi(x)
    \psi(y)|x-y|^{-a}\ud x\ud y,
  \]
  with the convention that $\|\varphi\|_{\cH_a} = |\int_{\bR}\varphi(x)\,\ud x|$
  (modulo the null space) if $a = 0$. Using the Cauchy-Schwarz inequality in
  $\cH_a$, and the fact that $\int_0^1 |t-s|^{-b}\ud s \le \frac{2}{1-b}$, it
  can be proved that for any function $\varphi \in \cH_{a,b}^{\otimes n}$
  supported in the time cube $[0,1]^n$ with finite right-hand side,
  \[
    \|\varphi\|_{\cH_{a,b}^{\otimes n}}^2 \le \left( \frac{2}{1-b}\right)^n
    \int_{[0,1]^n}\|\varphi(t_1,\ldots,t_n,\cdot)\|_{\cH_a^{\otimes n}}^2
    \ud t_1 \ldots \ud t_n.
  \]
  Indeed, after spatial Cauchy--Schwarz, bound the product of the two spatial
  norms by half the sum of their squares and integrate the time kernels.

  Using the scaling property \eqref{E:fnR-scale} of $f_{n,R}$ with $t = 1$ and
  $x = 0$, we see that
  \[
    \|f_{n,R}(t_1,\cdot,\ldots,t_n,\cdot,1,0)\|_{\cH_a^{\otimes
    n}}^2=R^{-an}\|f_{n}(t_1,\cdot,\ldots,t_n,\cdot,1,0)\|_{\cH_a^{\otimes
    n}}^2.
  \]
  Therefore,
  \begin{align}
    \nonumber
    n! \|\widetilde{f}_{n,R}(\cdot,1,0)\|_{\cH_{a,b}^{\otimes n}}^2 & \le
    \left( \frac{2}{1-b}\right)^n n!
    \int_{[0,1]^n}\|\widetilde{f}_{n,R}(t_1,\ldots,t_n,\cdot,1,0)\|_{\cH_a^{\otimes
    n}}^2 \ud t_1 \ldots \ud t_n \\
    \nonumber
    &= \left( \frac{2}{1-b}\right)^n R^{-an}
    \int_{\{0<t_1<\ldots<t_n<1\}}\|f_{n}(t_1,\ldots,t_n,\cdot,1,0)\|_{\cH_a^{\otimes
    n}}^2 \ud t_1 \ldots \ud t_n \\
    \nonumber
    &\le\left( \frac{2}{1-b}\right)^n R^{-an} c_a^n
    \int_{\{0<t_1<\ldots<t_n<1\}} \prod_{j=1}^{n}(t_{j+1}-t_j)^{2-a} \ud t_1
    \ldots \ud t_n \\
    \label{E:bound-fn}
    &=\left( \frac{2}{1-b}\right)^n R^{-an} c_a^n
    \frac{\Gamma(3-a)^n}{\Gamma\big((3-a)n+1\big)},
  \end{align}
  where $c_a > 0$ is a constant depending on $a$, and we denoted $t_{n+1} = 1$.
  For the spatial inequality, integrate the two earliest vertices first. With
  time gap $s$ and separation $d$ of the two following vertices, their
  contribution is
  \[
    \frac14\int_{-2s}^{2s}(2s-|y|)|d+y|^{-a}\,\ud y\le c_a s^{2-a}.
  \]
  The integral is maximal at $d = 0$: the two even decreasing factors have
  interval level sets whose overlaps are largest when centered together. Scaling
  at zero gives the bound, with $c_0 = 1$ when $a = 0$. Repeat for each
  successive pair. The conclusion follows from the chaos expansion of
  $M_R(cQ_{a,b})$ above and~\eqref{E:bound-fn}, letting $A_{a,b} = \frac{2}{1-b}
  c_a \Gamma(3-a)$ and using estimate \eqref{E:series-x}. This proves
  Lemma~\ref{L:MR-expo}.
\end{proof}

Now we are ready to prove Theorem~\ref{T:UB-th}.

\begin{proof}[Proof of Theorem~\ref{T:UB-th}]
  Let $T_N = Q-Q_N$. By the integral representation of the covariance, the
  quadratic form of $T_N$ on a real test $\varphi$ is $\int_{E\setminus
  E_N}|(T\varphi)(e)|^2\fm(\ud e)\ge0$. The kernel is also nonnegative since
  both kernel factors are nonnegative. By Lemma~\ref{L:log-subad},
  \begin{equation}\label{E:QT-bound}
    \log M_R(Q) \le \log M_R(Q_N)+\log M_R(T_N).
  \end{equation}

  We examine $T_N(t,x)$:
  \begin{align*}
    T_N(t,x)&= h_{\beta}(t)h_{\alpha}(x)- h_{\beta,N}(t)h_{\alpha,N}(x)\\
    &= \left(\int_{0}^{\frac{1}{N}} e^{-at^2} w_{\beta}(a)
    \ud a+\int_{N}^{\infty}
    e^{-at^2} w_{\beta}(a)\ud a \right)h_{\alpha}(x)\\
    &\quad +\left(\int_{0}^{\frac{1}{N}} e^{-ax^2} w_{\alpha}(a)
    \ud a+\int_{N}^{\infty} e^{-ax^2} w_{\alpha}(a)\ud a
    \right)h_{\beta,N}(t).
  \end{align*}

  Integrating $w_\rho$ up to $1/N$ and using $e^{-av^2}\le1$ gives the first
  estimate below. We use the following estimates: for any $v \in \bR$ and
  $\rho,\delta > 0$,
  \[
    \int_0^{\frac{1}{N}} e^{-av^2} w_{\rho}(a)\ud a \le
    \frac{N^{-\rho/2}}{\Gamma(1+\rho/2)} \quad \mbox{and} \quad
    \int_N^{\infty} e^{-av^2} w_{\rho}(a)\ud a \le N^{-\frac{\delta}{2}}
    \frac{\Gamma(\frac{\rho+\delta}{2})}{\Gamma(\frac{\rho}{2})}
    h_{\rho+\delta}(v).
  \]
  It follows that for any $\delta > 0$,
  \[
    T_N \le c \left(N^{-\frac{\beta}{2}} Q_{\alpha,0}+N^{-\frac{\delta}{2}}
    Q_{\alpha,\beta+\delta}+N^{-\frac{\alpha}{2}} Q_{0,\beta} +
    N^{-\frac{\delta}{2}} Q_{\alpha+\delta,\beta}\right),
  \]
  where $c > 0$ is a constant depending on $(\alpha,\beta,\delta)$. Note that
  $M_R(\rho_1) \le M_R(\rho_2)$ if $\rho_1 \le \rho_2$. Therefore, invoking also
  Lemma~\ref{L:log-subad}, we obtain
  \[
    \begin{aligned}
      \log M_R(T_N) &\le \log M_R\big(cN^{-\frac{\beta}{2}}
      Q_{\alpha,0}\big)+\log M_R\big(cN^{-\frac{\delta}{2}}
      Q_{\alpha,\beta+\delta}\big)\\ &\quad{}+\log
      M_R\big(cN^{-\frac{\alpha}{2}} Q_{0,\beta}\big) +\log
      M_R\big(c N^{-\frac{\delta}{2}} Q_{\alpha+\delta,\beta}\big).
    \end{aligned}
  \]

  Here $\delta$ is as chosen below, so Lemma~\ref{L:MR-expo} applies to all four
  terms. Using Lemma~\ref{L:MR-expo}, we obtain that
  \begin{equation}\label{E:T-bound}
    \frac{1}{R}\log M_R(T_N)\le \frac{\log c_1}{R}+c_2 e_N,
  \end{equation}
  where $e_N = (N^{-\frac{\beta}{2}})^{\frac{1}{3-\alpha}} +
  (N^{-\frac{\delta}{2}})^{\frac{1}{3-\alpha}} +
  (N^{-\frac{\alpha}{2}})^{\frac{1}{3}} +
  (N^{-\frac{\delta}{2}})^{\frac{1}{3-\alpha-\delta}}$, and $c_1,c_2 > 0$ are
  some constants depending on $(\alpha,\beta,\delta)$. We choose $0 < \delta <
  \min(1-\alpha,1-\beta)$ so that $e_N \to 0$ as $N \to \infty$.

  Using~\eqref{E:QT-bound} and~\eqref{E:T-bound}, followed by
  Theorem~\ref{T:UB-smooth}, we infer that
  \[
    \limsup_{R \to \infty} \frac{1}{R} \log M_R(Q) \le \limsup_{R \to \infty}
    \frac{1}{R} \log M_R(Q_N) + c_2 e_N \le 2D_N+c_2 e_N.
  \]
  Regard each $0 \le f \in \cK_N$ as zero off $E_N$. Cut off $y,z$ and smooth
  with a nonnegative mollifier to obtain $f_i \in \cF_+$ converging to $f$ in
  $\cK$: the scales may be kept in $[1/(2N),2N]$, where the weights are bounded
  above and below. By~\eqref{E:def-Cjk} on this band, $\|V_{f_i}-V_f^N\|_\infty
  \le C_{0,0,2N}\|f_i-f\|_{\cK}$. The Rayleigh suprema differ by at most this
  uniform error, so their square roots converge uniformly. Thus, the expression
  in the supremum at $f$ is the limit of values bounded by $D$. Take the
  supremum over $f$ and use $\log\|U_R^Q(1,0)\|_2 = \tfrac12\log M_R(Q)$ in the
  preceding display. The conclusion follows using the bound $D_N \le D$, and
  letting $N \to \infty$. This proves Theorem~\ref{T:UB-th}.
\end{proof}

\section{Upper bound for the \texorpdfstring{$p$}{p}-th moment}\label{S:upper-p}

In this section, we prove the upper bound in Theorem~\ref{T:main-th}.

\begin{theorem}\label{T:UB-th-p}
  Let $t_p$ be the value given in~\eqref{E:tp} and $D$ the constant
  given in~\eqref{E:def-D}. For any $p\ge2$, we have:
  \[
    \limsup_{t_p \to \infty} t_p^{-\frac{4-\alpha-\beta}{3-\alpha}} \log
    \|u_{\theta}(t,x)\|_p \le \theta^{\frac{1}{3-\alpha}}D.
  \]
\end{theorem}

\begin{proof}
  We use L\^e's moment comparison principle \eqref{E:hyper2}. Recalling
  relation~\eqref{E:u-rep1}, we obtain the following representation: for any $R
  > 0$
  \[
    \bE|u_{\theta}(t,x)|^2 =1+\sum_{n\geq 1} \big(\theta t^{4-\alpha-\beta}
    R^{\alpha}\big)^n n! \|\widetilde{f}_{n,R}(\cdot,1,0)\|_{\cH^{\otimes
    n}}^2.
  \]

  If we choose $R$ such that $\theta t^{4-\alpha-\beta} R^{\alpha} = R^3$, i.e.
  $R = \theta^{\frac{1}{3-\alpha}} t^{\frac{4-\alpha-\beta}{3-\alpha}}$, then
  the above expression coincides with $\bE|U_R^{Q}(1,0)|^2$
  (see~\eqref{E:2m-URQ}). It follows that for any $\theta > 0$,
  \[
    \bE|u_{\theta}(t,x)|^2=\bE\left|U_{\theta^{\frac{1}{3-\alpha}}
    t^{\frac{4-\alpha-\beta}{3-\alpha}}}^Q(1,0)\right|^2.
  \]
  We use this for $(p-1)\theta$ instead of $\theta$. Returning
  to~\eqref{E:hyper2}, it follows that
  \[
    \|u_{\theta}(t,x)\|_p
    \le \left\|u_{(p-1)\theta}(t,x)\right\|_2
    =   \left\|U_{((p-1)\theta)^{\frac{1}{3-\alpha}} t^{\frac{4-\alpha-\beta}{3-\alpha}}}^Q(1,0)\right\|_2
    =   \left\|U_{\theta^{\frac{1}{3-\alpha}} t_p^{\frac{4-\alpha-\beta}{3-\alpha}}}^Q(1,0)\right\|_2.
  \]
  Hence, by Theorem~\ref{T:UB-th},
  \[
    \theta^{-\frac{1}{3-\alpha}} \limsup_{t_p \to \infty}
    t_p^{-\frac{4-\alpha-\beta}{3-\alpha}} \log \|u_{\theta}(t,x)\|_p \le
    \limsup_{R\to \infty}\frac{1}{R} \log \|U_{R}^Q(1,0)\|_2 \le D.
  \]
  This proves Theorem~\ref{T:UB-th-p}.
\end{proof}

\section{Lower bound}\label{S:section-LB}

In this section, we prove the lower bound in Theorem~\ref{T:main-th}.

\medskip

We have the following preliminary result.

\begin{lemma}\label{L:V-HypA}
  For any $f \in \cF_{+}$, $V_{f}$ satisfies Hypothesis~\ref{H:hypA} (given in
  Section~\ref{S:section-wave}).
\end{lemma}

\begin{proof}
  Since $f$ has compact support in $E$, there exists an integer $N\geq2$ such
  that ${\rm supp}(f)\subset E_N$. Hence, $V_f=V_f^N$, and the conclusion
  follows from Lemma~\ref{L:VN-HypA}. This proves Lemma~\ref{L:V-HypA}.
\end{proof}

\begin{theorem}\label{T:LB-th-p}
  Let $D$ be the constant given in~\eqref{E:def-D}. For any $p \geq 2$,
  \[
    \liminf_{R \to \infty}\frac{1}{R} \log \|U_R^{(p-1)^{-1}Q}(1,0)\|_p \geq
    D.
  \]
  For every $\varepsilon > 0$ there is $R_0$ such that
  $R^{-1}\log\|U_R^{(p-1)^{-1}Q}(1,0)\|_p \ge D-\varepsilon$ for every $R \ge
  R_0$ and $p\ge2$. Consequently, by Lemma~\ref{L:u-U}, for the value
  $t_p$ given in~\eqref{E:tp}, we have:
  \[
    \liminf_{t_p \to \infty} t_p^{-\frac{4-\alpha-\beta}{3-\alpha}} \log
    \|u_{\theta}(t,x)\|_p \geq \theta^{\frac{1}{3-\alpha}}D.
  \]
\end{theorem}

\begin{proof}
  Recall from Lemma~\ref{L:T-map} that $W(\varphi) = B(T(\varphi))$ for any
  $\varphi \in \cH$, where $\{B(h);h \in \cK\}$ is a Gaussian white noise on
  $(E,\fm)$. Let $h \in \cK$ be arbitrary. Consider the random variable
  \begin{equation}\label{E:chaos-Zh}
    Z_h=e^{B(h)-\frac{1}{2}\|h\|_{\cK}^2}=\sum_{n\geq 0}\frac{1}{n!}
    I_n^{B}(h^{\otimes n}),
  \end{equation}
  where $I_n^B$ is the multiple integral of order $n$ with respect to $B$. Note
  that
  \begin{equation}\label{E:IWB}
    I_n^W(\varphi)=I_n^B(T^{\otimes n} \varphi) \quad \mbox{for all} \quad
    \varphi \in \cH^{\otimes n}.
  \end{equation}

  As in the proof of Proposition 6.3 of~\cite{balan.chen.ea:22:exact}, if
  $\frac{1}{p}+\frac{1}{q} = 1$, by H\"older's inequality,
  \[
    \bE[U_{R}^{Q}(1,0)Z_h] \le \|U_R^Q(1,0)\|_p \|Z_h\|_q=\|U_R^Q(1,0)\|_p
    \,
    e^{\frac{1}{2}(q-1)\|h\|_{\cK}^2},
  \]
  since $B(h) \sim N(0,\|h\|_{\cK}^2)$. We apply this for $Q$ replaced by
  $(p-1)^{-1}Q$. We obtain that
  \begin{equation}\label{E:LB-UpQ}
    \|U_R^{(p-1)^{-1}Q}(1,0)\|_p \geq \bE \left[ U_R^{(p-1)^{-1}Q}
    (1,0) Z_h\right] e^{-\frac{1}{2(p-1)}\|h\|_{\cK}^2}.
  \end{equation}

  We calculate the expectation on the right side. For this calculation take $h =
  \sqrt{R(p-1)}\,f$ with $f \in \cF_+$. The scale variables in these integrals
  lie in a compact interval. Since the integrands are nonnegative, Tonelli's
  theorem applies. Recalling the chaos expansion \eqref{E:chaos-URQp}, and
  using~\eqref{E:IWB}, it follows that:
  \[
    U_R^{(p-1)^{-1}Q} (1,0)=1+\sum_{n\geq 1} R^{\frac{3n}{2}}
    (p-1)^{-\frac{n}{2}} I_n^B\big(T^{\otimes n}
    f_{n,R}(\cdot,1,0)\big).
  \]
  Combining this with~\eqref{E:chaos-Zh}, and using orthogonality, we have:
  \begin{align*}
    \bE \left[ U_R^{(p-1)^{-1}Q} (1,0) Z_h\right]
    = 1+\sum_{n\geq 1} R^{\frac{3n}{2}} (p-1)^{-\frac{n}{2}}
    \langle
    T^{\otimes n} \widetilde{f}_{n,R}(\cdot,1,0),h^{\otimes n}
    \rangle_{\cK^{\otimes n}}.
  \end{align*}

  The integral formula for $T^{\otimes n}$ extends to the symmetrized kernel by
  the bounded compactly supported approximation used after~\eqref{E:fnR-scale}.
  The images converge in $\cK^{\otimes n}$ by isometry; they also converge
  pointwise by dominated convergence, since for each fixed tuple of features the
  product of their Gaussian kernels is bounded. Letting $\pmb{\xi} =
  (\xi_1,\ldots,\xi_n)$, we see that
  \begin{align*}
    & \langle T^{\otimes n} \widetilde{f}_{n,R}(\cdot,1,0),h^{\otimes n}
    \rangle_{\cK^{\otimes n}}\\
    &= \int_{E^n} \left(\int_{([0,1] \times
    \bR)^n} \widetilde{f}_{n,R} (\pmb{\xi},1,0)\prod_{j=1}^{n}K_{\xi_j}(e_j)
    \ud \pmb{\xi} \right) \prod_{j=1}^n h(e_j) \fm(\ud e_1)\ldots \fm(
    \ud e_n) \\
    &= \int_{([0,1] \times \bR)^n} \widetilde{f}_{n,R} (\pmb{\xi},1,0)
    \prod_{j=1}^{n} \left( \int_E K_{\xi_j} (e_j) h(e_j) \fm(\ud e_j)
    \right)
    \ud \pmb{\xi} \\
    &= \int_{([0,1] \times \bR)^n} \widetilde{f}_{n,R} (\pmb{\xi},1,0)
    \prod_{j=1}^{n} V_h(\xi_j) \ud \pmb{\xi}\\
    &=\int_{([0,1] \times
    \bR)^n} f_{n,R} (\pmb{\xi},1,0) \prod_{j=1}^{n} V_h(\xi_j)
    \ud \pmb{\xi},
  \end{align*}
  where for the last line we used the fact that $\langle \widetilde{f},g\rangle
  = \langle f, \widetilde{g} \rangle$. Hence,
  \begin{align*}
    \bE \left[ U_R^{(p-1)^{-1}Q} (1,0) Z_h\right]& =1+\sum_{n\geq 1} R^{2n}
    \int_{([0,1] \times \bR)^n} f_{n,R} (\pmb{\xi},1,0) \prod_{j=1}^{n}
    V_{\frac{h}{\sqrt{R(p-1)}}}(\xi_j) \ud \pmb{\xi} \\
    &= v_R \left( V_{\frac{h}{\sqrt{R(p-1)}}} ;1,0\right),
    \quad
    \mbox{for any $h=\sqrt{R(p-1)}\,f,\ f\in\cF_+$},
  \end{align*}
  or equivalently,
  \[
    \bE \left[ U_R^{(p-1)^{-1}Q} (1,0)
    Z_{\sqrt{R(p-1)}\,f}\right]=v_R(V_f;1,0) \quad \mbox{for any}
    \ f
    \in \cF_+.
  \]

  Coming back to~\eqref{E:LB-UpQ}, it follows that
  \[
    \|U_R^{(p-1)^{-1}Q}(1,0)\|_p \geq v_R(V_f;1,0) \,
    e^{-R\|f\|_{\cK}^2/2} \quad \mbox{for any} \ f \in
    \cF_+.
  \]
  Therefore, for any $f \in \cF_{+}$,
  \[
    \liminf_{R \to \infty}\frac{1}{R}\log \|U_R^{(p-1)^{-1}Q}(1,0)\|_p \geq
    \liminf_{R \to \infty} \frac{1}{R}
    \log v_R(V_f;1,0)-\frac{1}{2}\|f\|_{\cK}^2 \geq
    \Phi(V_f)-\frac{1}{2}\|f\|_{\cK}^2,
  \]
  where for the last inequality we used Theorem~\ref{T:lower-bound}, and
  Lemma~\ref{L:V-HypA}. The conclusion follows taking the supremum over all $f
  \in \cF_{+}$. At $p = 2$, the left side is bounded above by $\limsup_{R \to
  \infty}(2R)^{-1}\log M_R(Q) < \infty$ by Lemma~\ref{L:MR-expo} with $(a,b,c) =
  (\alpha,\beta,1)$. This bound is independent of $f$, so $D < \infty$. For
  strict positivity use the Riesz case of the second-moment lower bound cited in
  the introduction: Theorem~7.2 of~\cite{balan.conus:16:intermittency}, with
  $H=1-\beta/2$ and $a=\alpha$, applies since $0<\alpha<1$ implies $a<2$. With
  initial data $(1,0)$ and $\theta=H(2H-1)$ in our notation, it gives
  $\liminf_{t \to \infty}t^{-(4-\alpha-\beta)/(3-\alpha)}
  \log\|u_\theta(t,x)\|_2 > 0$. Theorem~\ref{T:UB-th-p} at $p = 2$ therefore
  gives $D > 0$.

  Given $\varepsilon > 0$, choose $f \in \cF_+$ for which the expression in the
  supremum exceeds $D-\varepsilon/2$. Before taking the liminf, the lower bound
  is $R^{-1}\log v_R(V_f;1,0)-\|f\|_{\cK}^2/2$, independent of $p$. By
  Theorem~\ref{T:log-v-th} it is at least $D-\varepsilon$ for every $R \ge
  R_0(f,\varepsilon)$, with the same threshold for all $p\ge2$. The upper bound
  through \eqref{E:hyper2} depends on $t,p$ only through $R_{p,t}$. Thus,
  \eqref{E:main-eq-tp} holds along every pair with $t_p \to \infty$. Since
  $\log\bE|u_\theta(t,x)|^p = p\log\|u_\theta(t,x)\|_p$, fixing $p$ gives
  \eqref{E:main-eq-t}; fixing $t$ gives \eqref{E:main-eq-p} because
  $p(p-1)^{1/(3-\alpha)}\sim p^{(4-\alpha)/(3-\alpha)}$. This proves
  Theorem~\ref{T:LB-th-p} and completes the proof of Theorem~\ref{T:main-th}.
\end{proof}

\appendix

\section{Cameron--Martin space}\label{S:appB}

In this section, we fix $N$ and study the Cameron--Martin space $\cH_R$ of the
Gaussian process $Z_R^* = \{Z_R^*(\xi);\ \xi \in C_R\}$ defined by
\[
  Z_R^*(\xi)=\frac{Z_N(\xi)}{\sqrt R},
  \qquad \xi\in C_R=[0,1]\times[-R,R].
\]
We regard $\mathbf{B} = C(C_R;\bC)$, equipped with the supremum norm, as a real
Banach space. Thus, $\mathbf{B}^*$ consists of the continuous real-linear
functionals from $\mathbf{B}$ to $\bR$. The process $Z_R^*$ is a centered
Gaussian random element in $\mathbf{B}$, that is, $\ell(Z_R^*)$ is a centered
real Gaussian random variable for every $\ell \in \mathbf{B}^*$.

We recall that $\cH_R$ is the completion of $\cH_{0,R} = \{S\ell;\ \ell \in
\mathbf{B}^*\}$ with respect to the inner product
\[
  \langle S\ell_1,S\ell_2\rangle_{\cH_R}
  =\bE[\ell_1(Z_R^*)\ell_2(Z_R^*)], \qquad \ell_1,\ell_2\in
  \mathbf{B}^*,
\]
where $S\ell = \bE[Z_R^*\ell(Z_R^*)]$ is a Bochner integral in $\mathbf{B}$.

We equip $\cK_N^{\bC} = \{f_1+if_2;\ f_1,f_2 \in \cK_N\}$ with the real inner
product
\[
  \langle f,g\rangle_{\cK_N^{\bC}}
  =\operatorname{Re}\int_{E_N} f(e)\overline{g(e)}\,\fm(\ud e).
\]
We extend \eqref{E:def-VfN} by setting $V_{f_1+if_2}^N = V_{f_1}^N+iV_{f_2}^N$.
For $f \in \cK_N^{\bC}$, let $T_Rf$ be the restriction of $V_f^N$ to $C_R$.
Since $\|K_\xi\|_{\cK_N}^2 = Q_N(0,0)$, Cauchy--Schwarz gives
\[
  \|T_Rf\|_\mathbf{B}\le \sqrt{Q_N(0,0)}\,\|f\|_{\cK_N^{\bC}}.
\]
Thus, $T_R:\cK_N^{\bC} \to \mathbf{B}$ is a bounded linear operator. By the
Riesz representation theorem, its adjoint $T_R^*:\mathbf{B}^* \to \cK_N^{\bC}$
satisfies
\begin{equation}\label{E:duality}
  \langle T_R^*\ell,f\rangle_{\cK_N^{\bC}}=\ell(T_Rf),
  \qquad \ell\in \mathbf{B}^*,\quad f\in\cK_N^{\bC}.
\end{equation}
Since $\ker(T_R)$ is closed, we have
\[
  \cK_N^{\bC}=\ker(T_R)\oplus\ker(T_R)^\perp \quad \text{and} \quad
  \overline{\operatorname{Range}(T_R^*)}=\ker(T_R)^\perp.
\]

\begin{lemma}
  \label{L:gaussian-series}We have the representation
  \[
    Z_R^*(\xi)=\frac{1}{\sqrt{2R}}\sum_{j\ge1}\gamma_j(T_RE_j)(\xi),
    \qquad \xi\in C_R,
  \]
  where $(\gamma_j)_{j\ge1}$ are independent standard real Gaussian random
  variables and $(E_j)_{j\ge1}$ is an orthonormal basis of the real Hilbert
  space $\cK_N^{\bC}$. The series converges in
  $L^2(\Omega;\mathbf{B})$.
\end{lemma}

\begin{proof}
  Let $(e_k)_{k\ge1}$ be an orthonormal basis of $\cK_N$. Recall that
  \[
    Z_R^*(\xi)=\frac{B_1(K_\xi)+iB_2(K_\xi)}{\sqrt{2R}},
  \]
  where $B_1$ and $B_2$ are independent real isonormal processes over
  $\cK_N$. By~\eqref{E:def-VfN},
  \[
    \langle K_\xi,e_k\rangle_{\cK_N}=V_{e_k}^N(\xi).
  \]
  Expanding $K_\xi$ in the basis $(e_k)_{k\ge1}$ and applying the
  isonormal isometry, we obtain, for each fixed $\xi \in C_R$,
  \[
    B_a(K_\xi)=\sum_{k\ge1}B_a(e_k)V_{e_k}^N(\xi)
    \quad\text{in }L^2(\Omega),\qquad a=1,2.
  \]

  For convergence in $L^2(\Omega;\mathbf{B})$, let $A$ be a bounded open
  rectangle containing $C_R$. By Parseval's identity and the derivative bounds
  for $K_\xi$, we have
  \[
    \begin{aligned}
      \sum_{k\ge1}\|V_{e_k}^N\|_{H^2(A)}^2
      &=\sum_{p+q\le2}\int_A\sum_{k\ge1}
      \left|\left\langle
      \partial_s^p\partial_x^qK_{s,x},e_k
      \right\rangle_{\cK_N}\right|^2\,\ud s\,\ud x\\
      &=\sum_{p+q\le2}\int_A
      \|\partial_s^p\partial_x^qK_{s,x}\|_{\cK_N}^2\,\ud s\,\ud x
      <\infty.
    \end{aligned}
  \]
  Since the $B_a(e_k)$ are independent standard Gaussians, for $M > L$,
  \[
    \bE\left\|
    \sum_{k=L+1}^{M}B_a(e_k)V_{e_k}^N
    \right\|_{H^2(A)}^2
    =\sum_{k=L+1}^{M}\|V_{e_k}^N\|_{H^2(A)}^2.
  \]
  By Sobolev embedding, there exists a constant $C_A > 0$ such that
  \[
    \bE\sup_{\xi\in C_R}
    \left|\sum_{k=L+1}^{M}B_a(e_k)V_{e_k}^N(\xi)\right|^2
    \le C_A\sum_{k>L}\|V_{e_k}^N\|_{H^2(A)}^2
    \longrightarrow0
    \qquad\text{as }L\to\infty.
  \]
  Thus, the partial sums converge in $L^2(\Omega;\mathbf{B})$. For each fixed
  $\xi \in C_R$, their limit equals $B_a(K_\xi)$ by uniqueness of the
  $L^2(\Omega)$ limit. Taking a countable dense subset of $C_R$ and using
  continuity, we identify this limit with the continuous version of
  $B_a(K_\cdot)$.

  Finally, set
  \[
    E_{2k-1}=e_k,\qquad E_{2k}=ie_k,\qquad
    \gamma_{2k-1}=B_1(e_k),\qquad \gamma_{2k}=B_2(e_k).
  \]
  Then $(E_j)_{j\ge1}$ is an orthonormal basis of $\cK_N^{\bC}$ over $\bR$, and
  the $\gamma_j$ are independent standard real Gaussian random variables. Since
  $T_R(ie_k) = iT_Re_k$, combining the two series gives the claimed
  representation in $L^2(\Omega;\mathbf{B})$. This proves
  Lemma~\ref{L:gaussian-series}.
\end{proof}

The following result describes the Cameron--Martin space $\cH_R$ and its norm.

\begin{theorem}\label{T:CM-th}
  \noindent a) $\cH_R = \operatorname{Range}(T_R)$;

  \smallskip
  \noindent b) For any $h \in \cH_R$,
  \[
    \|h\|_{\cH_R}^2
    =2R\inf_{f:\,T_Rf=h}\|f\|_{\cK_N^{\bC}}^2
    =2R\|f^*\|_{\cK_N^{\bC}}^2,
  \]
  where $f^*$ is the unique element of $\ker(T_R)^\perp$ satisfying $T_Rf^* =
  h$.
\end{theorem}

\begin{proof}
  \textit{a)} Let $\ell_1,\ell_2 \in \mathbf{B}^*$. Applying $\ell_1$ and
  $\ell_2$ to the series in Lemma~\ref{L:gaussian-series} and using
  $\bE[\gamma_j\gamma_k] = \delta_{jk}$, we get
  \begin{align*}
    \bE[\ell_1(Z_R^*)\ell_2(Z_R^*)]
    &=\frac{1}{2R}\sum_{j\ge1}\ell_1(T_RE_j)\ell_2(T_RE_j)\\
    &=\frac{1}{2R}\sum_{j\ge1}
    \langle T_R^*\ell_1,E_j\rangle_{\cK_N^{\bC}}
    \langle T_R^*\ell_2,E_j\rangle_{\cK_N^{\bC}}\\
    &=\frac{1}{2R}
    \langle T_R^*\ell_1,T_R^*\ell_2\rangle_{\cK_N^{\bC}},
  \end{align*}
  where we used \eqref{E:duality} and Parseval's identity. Hence, for any
  $\ell,m \in \mathbf{B}^*$,
  \[
    m(S\ell)=\bE[m(Z_R^*)\ell(Z_R^*)]
    =\frac{1}{2R}\langle T_R^*m,T_R^*\ell\rangle_{\cK_N^{\bC}}
    =m\!\left(\frac{1}{2R}T_RT_R^*\ell\right).
  \]
  Since $\mathbf{B}^*$ separates points of $\mathbf{B}$, this gives
  \[
    S\ell=\frac{1}{2R}T_RT_R^*\ell,
    \qquad
    \|S\ell\|_{\cH_R}^2
    =\bE[\ell(Z_R^*)^2]=\frac{1}{2R}\|T_R^*\ell\|_{\cK_N^{\bC}}^2.
  \]
  In particular, $\cH_{0,R} = T_R(\operatorname{Range}(T_R^*))$. For $g =
  T_R^*\ell$, we have $T_Rg = 2R\,S\ell$, and therefore
  \[
    \|T_Rg\|_{\cH_R}^2=(2R)^2\|S\ell\|_{\cH_R}^2
    =2R\|g\|_{\cK_N^{\bC}}^2.
  \]
  Thus, the map $T_R:\operatorname{Range}(T_R^*) \to \cH_{0,R}$ extends to a
  bijection from $\overline{\operatorname{Range}(T_R^*)}$ onto $\cH_R$. The
  inclusion $\cH_R\hookrightarrow\mathbf{B}$ is continuous: on the dense
  covariance range, Cauchy--Schwarz gives
  $\|S\ell\|_{\mathbf{B}}\le(\bE\|Z_R^*\|_{\mathbf{B}}^2)^{1/2}
  \|S\ell\|_{\cH_R}$, which extends by continuity. Since $T_R$ is bounded as a
  map into $\mathbf{B}$, this extension agrees with $T_R$. Hence,
  \[
    \cH_R
    =T_R\bigl(\overline{\operatorname{Range}(T_R^*)}\bigr)
    =T_R\bigl(\ker(T_R)^\perp\bigr).
  \]
  The orthogonal decomposition $\cK_N^{\bC} = \ker(T_R)\oplus\ker(T_R)^\perp$
  then gives $\cH_R = \operatorname{Range}(T_R)$. \smallskip

  \textit{b)} Given $h \in \cH_R$, let $f^* \in \ker(T_R)^\perp$ be the unique
  element satisfying $T_Rf^* = h$. Since $\overline{\operatorname{Range}(T_R^*)}
  = \ker(T_R)^\perp$, we may choose $f_n \in \operatorname{Range}(T_R^*)$ such
  that $f_n \to f^*$ in $\cK_N^{\bC}$. By the extension constructed in a),
  $T_Rf_n \to h$ in $\cH_R$. Applying the norm identity to $f_n$ and passing to
  the limit, we obtain
  \[
    \|h\|_{\cH_R}^2
    =\lim_{n\to\infty}\|T_Rf_n\|_{\cH_R}^2
    =2R\lim_{n\to\infty}\|f_n\|_{\cK_N^{\bC}}^2
    =2R\|f^*\|_{\cK_N^{\bC}}^2.
  \]

  If $T_Rf = h$, then $f-f^* \in \ker(T_R)$. By orthogonality,
  \[
    \|f\|_{\cK_N^{\bC}}^2
    =\|f^*\|_{\cK_N^{\bC}}^2+\|f-f^*\|_{\cK_N^{\bC}}^2
    \ge\|f^*\|_{\cK_N^{\bC}}^2.
  \]
  Equality holds exactly when $f = f^*$. Therefore,
  \[
    \|h\|_{\cH_R}^2
    =2R\inf_{f:\,T_Rf=h}\|f\|_{\cK_N^{\bC}}^2
    =2R\|f^*\|_{\cK_N^{\bC}}^2.
  \]
  This proves Theorem~\ref{T:CM-th}.
\end{proof}

\begin{corollary}\label{C:CM-ball}
  Let $B_{\cH_R}$ be the open unit ball in $\cH_R$. For any $a > 0$, the open
  ball of radius $\sqrt{2R}a$ in $\cH_R$ admits the representation:
  \[
    \sqrt{2R} a B_{\cH_R} = \{T_Rf; f \in
    \cK_{N}^{\bC},\|f\|_{\cK_N^{\bC}}<a\}.
  \]
\end{corollary}

\begin{proof}
  First, let $h \in \cH_{R}$ be such that $\|h\|_{\cH_R} < \sqrt{2R} a$. By
  Theorem~\ref{T:CM-th}.b), there exists $f^* \in \cK_{N}^{\bC}$ with $T_R f^* =
  h$ such that $\|h\|_{\cH_R} = \sqrt{2R} \|f^*\|_{\cK_N^{\bC}}$. Hence,
  $\|f^*\|_{\cK_N^{\bC}} < a$.

  Secondly, let $h = T_Rf$ for some $f \in \cK_{N}^{\bC}$ with
  $\|f\|_{\cK_N^{\bC}} < a$. By Theorem~\ref{T:CM-th}.b),
  \[
    \|h\|_{\cH_R}^2 =2R \inf_{g :T_R g=h }\|g\|_{\cK_N^{\bC}}^2 \le 2R
    \|f\|_{\cK_N^{\bC}}^2 < 2R a^2,
  \]
  and hence $\|h\|_{\cH_R} < \sqrt{2R}a$, i.e. $h \in \sqrt{2R} a B_{\cH_R}$.
  This proves Corollary~\ref{C:CM-ball}.
\end{proof}

\section{Comparison with the time-independent case}\label{S:static-comparison}

In this section, we show that when $\beta=0$, the limit in Theorem \ref{T:main-th} coincides with the one obtained in \cite{balan.chen.ea:22:exact} in the time-independent (or static) case, when $d=1$ and $\gamma(x)=|x|^{-\alpha}$.

\medskip

Let $\cK_{\mathrm{stat}} \coloneqq L^2\big((0,\infty)\times\bR,\,
w_{\alpha}(b)\ud b\ud z\big)$. For any nonnegative $f \in
C_c^{\infty}((0,\infty)\times\bR)$, we let
\[
  V_f^{\mathrm{stat}}(x) \coloneqq \int_0^\infty\int_{\bR}
  \fg_b(x-z)f(b,z)w_{\alpha}(b)\ud z\ud b, \quad x \in \bR.
\]
Comparing this with the definition of $V_f(t,x)$, we see that the temporal factor
$\fg_a(t-y)$ in $K_{t,x}$ is replaced by 1.
We define the time-independent (or static) variational constant:
\begin{equation}\label{E:D-stat}
  D_{\mathrm{stat}} \coloneqq \sup_{f}
  \left\{\Phi\big(V_f^{\mathrm{stat}}\big)
  -\frac{1}{2}\|f\|_{\cK_{\mathrm{stat}}}^2\right\},
\end{equation}
where the supremum is over all nonnegative $f \in
C_c^{\infty}((0,\infty)\times\bR)$, and $V_f^{\mathrm{stat}}$ is viewed as a potential on $[0,1]\times\bR$ that does not depend on time.

\begin{proposition}
  \label{P:static-comparison}
  We have:
  \begin{equation}\label{E:D-stat-explicit}
    D_{\mathrm{stat}} = 2^{-\frac{\alpha}{2(3-\alpha)}}\, \frac{3-\alpha}{2}
    \left(\frac{2\sqrt{\cM}}{4-\alpha}\right)^{\frac{4-\alpha}{3-\alpha}},
  \end{equation}
  where $\cM$ is the variational constant $\cM(\gamma,1)$
  of~\cite{balan.chen.ea:22:exact}:
  \begin{equation}\label{E:cM}
    \cM \coloneqq \sup_{g\in H^1(\bR),\|g\|_{L^2(\bR)}=1}
    \left\{\left(\int_{\bR}\int_{\bR}
    |x-y|^{-\alpha}g(x)^2g(y)^2\ud x\ud y\right)^{\frac{1}{2}}
    -\frac{1}{2}\|g'\|_{L^2(\bR)}^2\right\}.
  \end{equation}
  In particular, $\theta^{1/(3-\alpha)}D_{\mathrm{stat}}$ equals the constant in
  Theorem~1.6(i) of~\cite{balan.chen.ea:22:exact}. The time exponent and
  value $t_p$ in Theorem~\ref{T:main-th} also agree with those of the
  static model upon setting $\beta=0$.
\end{proposition}

\begin{proof}
  Theorem~\ref{T:main-th} is the time-dependent counterpart of Theorem 1.6
  of~\cite{balan.chen.ea:22:exact} (with $d = 1$). When $\beta = 0$, the
  temporal covariance is $1$, and the exponent and $t_p$ become
  $(4-\alpha)/(3-\alpha)$ and $(p-1)^{1/(4-\alpha)}t$, respectively, as
  in~\cite{balan.chen.ea:22:exact}. Since the definition of $w_\beta$ requires
  $\beta > 0$, we compare the constants by considering the static variational
  problem separately.

  For this comparison, let $f$ range over the nonnegative functions in
  $C_c^\infty((0,\infty)\times\bR)$. To lighten the notation, in this proof we
  write $V_f$ for $V_f^{\mathrm{stat}}$ and $\|f\|$ for
  $\|f\|_{\cK_{\mathrm{stat}}}$, both defined at the beginning of this
  section; that is,
  \[
    V_f(x)=\int_0^\infty\int_{\bR}
    \fg_b(x-z)f(b,z)w_\alpha(b)\,\ud z\,\ud b,
    \qquad
    \|f\|^2=\int_0^\infty\int_{\bR}
    |f(b,z)|^2w_\alpha(b)\,\ud z\,\ud b.
  \]
  Since $V_f$ does not depend on time, we write $\lambda(V_f) =
  \lambda_{V_f}(t)$ and have $\Phi(V_f) = \sqrt{\lambda(V_f)}$. Thus, the
  definition~\eqref{E:D-stat} of the static variational constant reads
  \[
    D_{\mathrm{stat}}
    =\sup_{f\ge0}
    \left\{\sqrt{\lambda(V_f)}-\frac12\|f\|^2\right\}.
  \]

  For $g \in H^1(\bR)$ with $\|g\|_2 = 1$, set
  \[
    A(g)^2=\int_{\bR}\int_{\bR}
    |x-y|^{-\alpha}g(x)^2g(y)^2\,\ud x\,\ud y,
    \qquad B(g)=\|g'\|_2^2.
  \]
  With this notation, the constant $\cM$ in~\eqref{E:cM} is
  \[
    \cM=\sup_{\substack{g\in H^1(\bR),\;\|g\|_2=1}}
    \left\{A(g)-\frac12B(g)\right\}.
  \]

  Fix $r>0$ and restrict to controls with $\|f\|\le r$. For a fixed $g$, the
  function $(b,z)\mapsto\int_{\bR}\fg_b(x-z)g(x)^2\,\ud x$ has norm $A(g)$ in
  $\cK_{\mathrm{stat}} $, by Fubini's theorem and
  the spatial part of~\eqref{E:KK-int}, namely,
  $\int_0^\infty\int_{\bR}\fg_b(x-z)\fg_b(y-z)w_\alpha(b)\,\ud z\,\ud b =
  h_\alpha(x-y) = |x-y|^{-\alpha}$ for $x \ne y$ (see~\eqref{E:semigroup}). By
  Fubini's theorem and Cauchy--Schwarz,
  \begin{align*}
    \int_{\bR}V_f(x)g(x)^2\,\ud x
    &=\int_{\bR}
    \left(\int_0^\infty\int_{\bR}
    \fg_b(x-z)f(b,z)w_\alpha(b)\,\ud z\,\ud b\right)
    g(x)^2\,\ud x\\
    &=\int_0^\infty\int_{\bR}
    f(b,z)\left(\int_{\bR}\fg_b(x-z)g(x)^2\,\ud x\right)
    w_\alpha(b)\,\ud z\,\ud b\\
    &=\left\langle f,
    \int_{\bR}\fg_b(x-z)g(x)^2\,\ud x\right\rangle\\
    &\le\|f\|
    \left\|\int_{\bR}\fg_b(x-z)g(x)^2\,\ud x\right\|
    =\|f\|A(g)\le rA(g).
  \end{align*}
  By approximation in this norm by nonnegative smooth compactly supported
  functions, the bound is sharp. Therefore,
  \begin{equation}\label{E:sup-f-stat}
    \sup_{\substack{f\ge0,\;\|f\|\le r}}
    \int_{\bR}V_f(x)g(x)^2\,\ud x=rA(g).
  \end{equation}
  For $s > 0$, let $g_s(x) = s^{1/2}g(sx)$. A change of variables gives
  \[
    \|g_s\|_2=1,\qquad
    A(g_s)=s^{\alpha/2}A(g),\qquad
    B(g_s)=s^2B(g).
  \]
  Taking $s = (r/2)^{2/(4-\alpha)}$, so that $rs^{\alpha/2} = 2s^2$, yields
  \[
    rA(g_s)-B(g_s)
    =2^{-\alpha/(4-\alpha)}r^{4/(4-\alpha)}
    \left(A(g)-\frac12B(g)\right).
  \]
  For each fixed $s > 0$, the map $g\mapsto g_s$ is a bijection of $\{g \in
  H^1(\bR):\|g\|_2 = 1\}$ onto itself. Taking the supremum over $g_s$ is the
  same as taking it over $g$. It follows that
  \[
    \sup_{\substack{f\ge0,\; \|f\|\le r}}\lambda(V_f)
    =2^{-\alpha/(4-\alpha)}r^{4/(4-\alpha)}\cM.
  \]

  For each fixed $f$, we have
  \[
    \sqrt{\lambda(V_f)}-\frac12\|f\|^2
    =\sup_{r\ge\|f\|}
    \left\{\sqrt{\lambda(V_f)}-\frac{r^2}{2}\right\}.
  \]
  Taking the supremum over $f$ and exchanging the two suprema gives
  \[
    \begin{aligned}
      D_{\mathrm{stat}}
      &=\sup_{r\ge0}
      \left\{
      \sqrt{\sup_{\substack{f\ge0\\\|f\|\le r}}\lambda(V_f)}
      -\frac{r^2}{2}
      \right\}\\
      &=\sup_{r\ge0}
      \left\{
      2^{-\alpha/[2(4-\alpha)]}\sqrt{\cM}\,
      r^{2/(4-\alpha)}-\frac{r^2}{2}
      \right\}\\
      &=2^{-\alpha/[2(3-\alpha)]}\frac{3-\alpha}{2}
      \left(\frac{2\sqrt{\cM}}{4-\alpha}\right)^{(4-\alpha)/(3-\alpha)}.
    \end{aligned}
  \]
  The last equality follows by differentiating with respect to $r$ and setting
  the derivative equal to zero. This proves~\eqref{E:D-stat-explicit}, and
  hence Proposition~\ref{P:static-comparison}.
\end{proof}

\begin{remark}\label{R:lambda-M}
  The quantities $\lambda_V(t)$ in~\eqref{E:lambda-V} and $\cM$
  in~\eqref{E:cM} are of the same type: each is the
  supremum, over $g\in H^1(\bR)$ with $\|g\|_{L^2(\bR)}=1$, of a potential
  energy minus a kinetic energy. In $\lambda_V(t)$ the potential $V(t,\cdot)$ is
  given. For fixed $g$, optimizing the potential energy over the nonnegative
  test functions $f$ in~\eqref{E:D-stat} gives
  \[
    \sup_{f\ge 0,\: \|f\|_{\cK_{\mathrm{stat}}}\le 1}
    \int_{\bR}V_f^{\mathrm{stat}}(x)g(x)^2\ud x
    =\left(\int_{\bR}\int_{\bR}
    |x-y|^{-\alpha}g(x)^2g(y)^2\ud x\ud y\right)^{\frac{1}{2}},
  \]
  which is the first term in~\eqref{E:cM}; see~\eqref{E:sup-f-stat} with $r =
  1$. For each control radius $r$, the proof of
  Proposition~\ref{P:static-comparison} first optimizes over $f$ for fixed $g$ and
  rescales $g$; it then optimizes over $r$ to obtain~\eqref{E:D-stat-explicit}. In the time-dependent case, the same
  control $f$ enters $\lambda_{V_f}(t)$ at all times $t\in[0,1]$, and we do not
  know of a similar reduction for $D$ in~\eqref{E:def-D}.
\end{remark}

\section*{Acknowledgments}
\addcontentsline{toc}{section}{Acknowledgments} Jinxin Wang thanks
the Department of Mathematics and Statistics at Auburn University for its
hospitality during his research visit. Le Chen thanks OpenAI for providing
access to its AI tools through the ChatGPT for Academic Researchers program.

\paragraph{Funding.}
Raluca M.~Balan was supported by a grant from the Natural
Sciences and Engineering Research Council of Canada. Le Chen was partially
supported by the National Science Foundation through CAREER award
\mbox{DMS-2443823} and by the Simons Foundation through award
\mbox{MPS-TSM-00959981}. A Mitacs Globalink Research Award Abroad (IT52612)
supported Jinxin Wang's visit to Auburn University.

\medskip

{\bf Statement on AI use.} The proof strategy and a
first version of the arguments were obtained by the authors with the help
of an AI system. The validity of all proofs was checked by the authors, who take full
responsibility for the content. Arriving at the argument required many exchanges, in which the authors examined the proposed steps, and asked follow-up questions. Even after the main idea emerged, the argument was difficult to follow in the form that was provided. The authors spent a substantial amount of time and effort reconstructing the proof, checking its details, and deciding how to explain it. The final presentation is the illustration of that work.

\end{document}